%% file: main.tex
\documentclass[smallextended,referee,envcountsect,envcntsame]{svjour3}
\smartqed
\usepackage{amsmath,amssymb,amsfonts,stmaryrd,latexsym}
\usepackage{graphicx}
\usepackage{mathtools}
\usepackage{bm}
\usepackage{autonum}
\usepackage{algpseudocode}
\usepackage{algpseudocode}
\usepackage{algorithm}
\usepackage{enumitem}

\newcommand{\R}{\mathbb{R}}
\newcommand{\N}{\mathbb{N}}

\newcommand{\argmin}[1][]{\underset{{#1}}{\mathrm{argmin}}}
\newcommand{\product}[3][]{\left\langle {#2} , {#3} \right\rangle_{#1}}

\newcommand{\vv}{\bm{v}}
\newcommand{\w}{\bm{w}}
\newcommand{\x}{\bm{x}}
\newcommand{\y}{\bm{y}}

\newcommand{\f}{\bm{f}}
\newcommand{\bxi}{\bm{\xi}}
\newcommand{\bzeta}{\bm{\zeta}}
\newcommand{\Li}{\bm{L}}

\newcommand{\G}{\bm{G}}

\newcommand{\A}{\bm{A}}
\newcommand{\B}{\bm{B}}
\newcommand{\C}{\bm{C}}

\newcommand{\Q}{\bm{Q}}
\newcommand{\TFBF}{\bm{T}_{\mathrm{FBF}}}
\newcommand{\Talpha}{\bm{T}^{\alpha}_{\FBF}}
\newcommand{\uf}[1]{\bm{\mathfrak{f}}_{{#1}}^{\langle \mathrm{u} \rangle}}
\newcommand{\uG}{\bm{\mathfrak{G}}}
\newcommand{\st}{\mathrm{s.t.} \ }
\newcommand{\find}{\mathrm{find} \ }
\newcommand{\Id}{\mathrm{Id}}
\newcommand{\FBF}{\mathrm{FBF}}
\newcommand{\FB}{\mathrm{FB}}
\newcommand{\Fix}{\mathrm{Fix}}
\newcommand{\VI}{\mathrm{VI}}
\newcommand{\zer}{\mathrm{zer}}

\newcommand{\dom}{\mathrm{dom}}
\newcommand{\uV}{\bm{\mathcal{V}}^{\langle \mathrm{u} \rangle}}

\newcommand{\T}{\bm{T}}

\newcommand{\clBall}{\overline{B}(0;r)}
\newcommand{\SGNE}[6]{\bm{\mathcal{S}}_{\mathrm{GNE}}^{#1, #2}((#3, #4)_{i\in\mathcal{I}}, #5, #6)}
\newcommand{\sSGNE}[2]{\bm{\mathcal{S}}_{\mathrm{GNE}}^{#1, #2}}

\journalname{}

\begin{document}

\title{Hierarchical Variational Inequality Problem for Noncooperative Game-Theoretic Selection of \\Generalized Nash Equilibrium}

\author{Shota Matsuo \and Keita Kume \and \\Isao Yamada}

\institute{
   Shota Matsuo 
   \at Department of Information and Communications Engineering, Institute of Science Tokyo, S3-60, Ookayama, Meguro-ku, Tokyo 152-8552, Japan, e-mail: matsuo@sp.ict.e.titech.ac.jp\\ \\
   Keita Kume 
   \at Department of Information and Communications Engineering, Institute of Science Tokyo, S3-60, Ookayama, Meguro-ku, Tokyo 152-8552, Japan, e-mail: kume@sp.ict.e.titech.ac.jp\\ \\
   Isao Yamada 
   \at Department of Information and Communications Engineering, Institute of Science Tokyo, S3-60, Ookayama, Meguro-ku, Tokyo 152-8552, Japan, e-mail: isao@sp.ict.e.titech.ac.jp
}

\date{Received: date / Accepted: date}

\maketitle

\begin{abstract}
   The equilibrium selection problem in the variational Generalized Nash Equilibrium Problem (v-GNEP) has been reported as an optimization problem defined over the solution set of v-GNEP, called in this paper the lower-level v-GNEP. 
   However, to make such a selection fair for every player, we have to rely on an unrealistic assumption, that is, the availability of a trusted center that does not induce any bias among players. 
   In this paper, to ensure fairness for every player even in the process of equilibrium selection, we propose a new equilibrium selection problem, named the upper-level v-GNEP. 
   The proposed upper-level v-GNEP is formulated as a v-GNEP defined over the solution set of the lower-level v-GNEP. 
   We also present an iterative algorithm, of guaranteed convergence to a solution of the upper-level v-GNEP, as an application of the hybrid steepest descent method to a fixed point set characterization of the solution of the lower-level v-GNEP.  
   Numerical experiments illustrate the proposed equilibrium selection and algorithm. 
\end{abstract}

\input{docs/Introduction.tex}
\input{docs/Preliminaries.tex}
\input{docs/Algorithm.tex}
\input{docs/Application.tex}

\section{Conclusions}
\label{sec:conclusion}
We proposed a new formulation of equilibrium selection problem, named the upper-level v-GNEP, which was designed for a fair equilibrium selection without assuming any trusted center. 
We also proposed an iterative algorithm for the upper-level v-GNEP as an application of the hybrid steepest descent method to a fixed point set characterization of the solution of the lower-level v-GNEP. 
Numerical experiments illustrate the proposed equilibrium selection and algorithm. 

\begin{acknowledgements}\.
\\ This work was supported by JSPS Grants-in-Aid (19H04134, 24K23885). 
\end{acknowledgements}

\appendix  

\input{docs/Appendix.tex}

\bibliographystyle{spmpsci}
\bibliography{ref}

\end{document}

%% file: docs/Introduction.tex
\section{Introduction}
Game theory dates back to the pioneering work of von Neumann and Morgenstern \cite{NM1944}, and the innovative idea of \emph{Nash equilibrium (NE)}, introduced by John Nash \cite{nash1950equilibrium,nash1951non}, triggered the drastic expansion of applications of game theory. 
The NE is a well-balanced solution of non-cooperative games, in which multiple decision-makers $\mathcal{I}\coloneqq\{1, \dots, m\}$, called players, aim to decrease respectively their cost functions as much as possible. 
At the NE, any player $i\in\mathcal{I}$ cannot decrease solely $i$'s cost function by changing $i$'s variable, called $i$'s strategy, as long as the other players' strategies are unchanged. 
The NE has been generalized \cite{facchineiGeneralizedNashReview2010,buiNashEquilibrium2021} and advanced \cite{luo1996mathematical,g.scutariEquilibriumSelectionPower2012,g.scutariRealComplexMonotone2014,e.benenatiOptimalSelectionTracking2023,w.heDistributedOptimalVariational2024a} toward one of ideal goals of decision processes in a variety of modern engineering/social systems (see, e.g., \cite{yeDistributedNashEquilibrium2023}), such as wireless communication \cite{a.deligiannisGameTheoreticPowerAllocation2017,w.wangNonConvexGeneralizedNash2022,y.chenQoEAwareDecentralizedTask2024}, smart grids \cite{i.atzeniNoncooperativeDayAheadBidding2014,zhengPeertopeerEnergyTrading2022}, and machine learning \cite{NIPS2014_5ca3e9b1,h.tembineDeepLearningMeets2020}. 
In this paper, we start with the following \emph{Generalized Nash Equilibrium Problem (GNEP)} as an extension of the classical \emph{Nash equilibrium problem (NEP)} (see Remark \ref{rem:GNEP_NEP}). 
\begin{problem}[{Generalized Nash Equilibrium Problem (GNEP), see, e.g., \cite[Def. 3.6]{facchineiGeneralizedNashReview2010}, \cite[Exm. 3.4]{buiNashEquilibrium2021}}]
  \label{prob:GNEP}
  Consider a non-cooperative game among players in $\mathcal{I}\coloneqq\{1, \dots, m\}$ with $m\in\N\setminus\{0\}$. 
  We follow the notations used in \cite{buiNashEquilibrium2021} as:
  \begin{enumerate}[label=\textup{(N\arabic*)}, leftmargin=2.5em]
    \item 
    For every $i \in \mathcal{I}$, the strategy of player $i$ is defined by $x_i \in \mathcal{H}_i$, where $(\mathcal{H}_i, \product[\mathcal{H}_i]{\cdot}{\cdot}, \|\cdot\|_{\mathcal{H}_i})$ is a finite-dimensional real Hilbert space. 
    \item
    Whole players' strategies\footnote{Multiple Hilbert spaces $(\mathcal{H}_i, \langle \cdot, \cdot \rangle_{\mathcal{H}_i}, \|\cdot\|_{\mathcal{H}_i})$ $(i\in\mathcal{I})$ can be used to build a new Hilbert space $\bm{\mathcal{H}}\coloneqq\bigl(\textstyle\bigtimes_{i\in\mathcal{I}}\mathcal{H}_i\bigr)\coloneqq \mathcal{H}_1 \times \mathcal{H}_2 \times \cdots \times \mathcal{H}_m = \bigl\{ \x \coloneqq (x_i)_{i\in\mathcal{I}} \coloneqq (x_1, x_2, \ldots, x_m )\,\bigl|\, x_i \in \mathcal{H}_i\,(i\in\mathcal{I})\bigr\}$ equipped with (i) the addition $\bm{\mathcal{H}}\times\bm{\mathcal{H}}\to\bm{\mathcal{H}}: \bigl( (x_i)_{i\in\mathcal{I}},(y_i)_{i\in\mathcal{I}}\bigr) \mapsto (x_i + y_i)_{i\in\mathcal{I}}$, (ii) the scalar multiplication $\R\times\bm{\mathcal{H}}\to\bm{\mathcal{H}}:\bigl( \alpha, (x_i)_{i\in\mathcal{I}}\bigr) \mapsto (\alpha x_i)_{i\in\mathcal{I}}$, and (iii) the inner product $\bm{\mathcal{H}}\times\bm{\mathcal{H}}\to\mathbb{R}:\bigl( (x_i)_{i\in\mathcal{I}}, (y_i)_{i\in\mathcal{I}} \bigr) \mapsto \bigl\langle (x_i)_{i\in\mathcal{I}}, (y_i)_{i\in\mathcal{I}} \bigr\rangle_{\bm{\mathcal{H}}} \coloneqq \sum_{i\in\mathcal{I}} \langle x_i, y_i \rangle_{\mathcal{H}_i}$ and its induced norm $\bm{\mathcal{H}}\to\mathbb{R}:\x\mapsto\|\x\|_{\bm{\mathcal{H}}}\coloneqq\sqrt{\langle \x,\x \rangle_{\bm{\mathcal{H}}}}$.} and players' strategies other than player $i\in \mathcal{I}$ are denoted respectively by $\x\coloneqq(x_1, \dots, x_m) \in \bm{\mathcal{H}} \coloneqq \textstyle \bigtimes_{i\in\mathcal{I}}\mathcal{H}_i$ and by 
    $\x_{\smallsetminus i}\coloneqq(x_1, \dots, x_{i-1}, x_{i+1}, \dots, x_m)$. 
    \item
    For each $i \in \mathcal{I}$ and any $(x_i, \y) \in \mathcal{H}_i \times \bm{\mathcal{H}}$, $(x_i;\y_{\smallsetminus i})$ stands for \linebreak $(y_1, \dots, y_{i-1}, x_i, y_{i+1}, \dots, y_m) \in \bm{\mathcal{H}}$. 
    \item 
    $(\mathcal{G}, \product[\mathcal{G}]{\cdot}{\cdot}, \|\cdot\|_{\mathcal{G}})$ is a finite-dimensional real Hilbert space. 
  \end{enumerate}
  Then, the \emph{Generalized Nash Equilibrium Problem (GNEP)} is formulated as: 
  \begin{align}
    \label{eq:GNEP}
    &\mathrm{find} \ {\x} \in \bm{\mathcal{S}}_{\mathrm{GNE}}^{\bm{\mathcal{H}}, \mathcal{G}}((\f_i, C_i)_{i\in\mathcal{I}}, \Li, D)\coloneqq\Bigl\{(\overline{x}_1, \dots, \overline{x}_m)\in \bm{\mathcal{H}} \Bigm\vert \\
    &\qquad \qquad \qquad \qquad (\forall i \in I) \ \overline{x}_i \in \argmin[x_i \in C_i] \ \f_i(x_i; \overline{\x}_{\smallsetminus i}) \ \st \Li (x_i ;\overline{\x}_{\smallsetminus i}) \in D\Bigr\}, 
  \end{align}
  where 
  \begin{enumerate}[label=\textup{(\alph*)}]
    \item \label{item:Li_C_D} $\Li:\bm{\mathcal{H}}\to\mathcal{G}$ is a linear operator, and  $C_i\subset\mathcal{H}_i \ (i\in\mathcal{I})$ and $D\subset\mathcal{G}$ are closed convex sets satisfying 
    \begin{align}
      \label{eq:fracture_C}
      \bm{\mathfrak{C}} \coloneqq \C\cap\Li^{-1}(D) \coloneqq \bigl\{\x \in \textstyle \bigtimes_{i \in \mathcal{I}} C_i \bigm\vert \Li\x \in D\bigr\} \neq \varnothing, 
    \end{align}
    where 
    $\C\coloneqq\bigtimes_{i \in \mathcal{I}} C_i$ and  $\Li^{-1}(D)\coloneqq\{\x\in\bm{\mathcal{H}} \vert \Li\x\in D\}$. 
    \item \label{item:cost_low} For every $i\in\mathcal{I}$, $\f_i:\bm{\mathcal{H}}\to\R$ satisfies that, for every $\x \in \bm{\mathcal{H}}$, $\f_i(\cdot;\x_{\smallsetminus i}):\mathcal{H}_i\to\R$ is convex and differentiable over $\mathcal{H}_i$. 
  \end{enumerate}
\end{problem}

\begin{remark}
  \label{rem:GNEP_NEP}
  Consider a simple case of Problem \ref{prob:GNEP} where $(\mathcal{G}, \Li)\coloneqq(\bm{\mathcal{H}}, \Id)$ and $D\coloneqq\bigtimes_{i\in\mathcal{I}} D_i$ with nonempty closed convex sets $D_i\subset\mathcal{H}_i \ (i\in\mathcal{I})$. 
  In this case, Problem \ref{prob:GNEP} is the so-called Nash Equilibrium Problem (NEP) (see, e.g., \cite[Sec. 1.4.2]{FPvariational2003}). 
\end{remark}

Facchinei, Fischer, and Piccialli \cite{facchineiGeneralizedNashGames2007} introduced a variational inequality $\mathrm{VI}(\bm{\mathfrak{C}}, \G)$:
\begin{align}
  \label{eq:VE}
  \find \x \in \bm{\mathcal{V}}\coloneqq\left\{\vv \in \bm{\mathfrak{C}} \mid (\forall \w \in \bm{\mathfrak{C}}) \ \product[\bm{\mathcal{H}}]{\G(\vv)}{\w-\vv} \geq 0 \right\} \subset\bm{\mathcal{H}}, 
\end{align}
where $\G$ is defined, with gradients $\nabla_i \f_i(\cdot;\x_{\smallsetminus i}):\mathcal{H}_i\rightarrow\mathcal{H}_i$ of $\f_i(\cdot;\x_{\smallsetminus i})$, as
\begin{equation}
  \label{eq:G}
  \bm{G}:\bm{\mathcal{H}}\rightarrow\bm{\mathcal{H}}:\x\mapsto(\nabla_1 \f_1(\x), \cdots, \nabla_m \f_m(\x)).
\end{equation}

Indeed, the solution set $\bm{\mathcal{V}}$ of $\mathrm{VI}(\bm{\mathfrak{C}}, \G)$ is a special subset of \linebreak $\SGNE{\bm{\mathcal{H}}}{\mathcal{G}}{\f_i}{C_i}{\Li}{D}$ (see Fact \ref{fact:vGNE_GNEP}(ii)). 
A point in $\bm{\mathcal{V}}$ of \eqref{eq:VE} has been revealed to possess more desirable properties, such as \emph{fairness} and \emph{larger social stability}, than $\sSGNE{\bm{\mathcal{H}}}{\mathcal{G}}\setminus\bm{\mathcal{V}}$ \cite{facchineiGeneralizedNashReview2010,kulkarniVariationalEquilibriumRefinement2012}, and referred to as a \emph{variational equilibrium} \cite[Def. 3.10]{facchineiGeneralizedNashReview2010} or \emph{variational GNE (v-GNE)} \cite{BelgioiosoSemiDecentralizedNash2017}. 
In this paper, problem \eqref{eq:VE} is specially referred to as the \underline{\emph{v-GNE Problem (v-GNEP)} for $\sSGNE{\bm{\mathcal{H}}}{\mathcal{G}}$ introduced in \eqref{eq:GNEP}}. 
Recent applications of v-GNE are found, e.g., in distributed control \cite{g.belgioiosoDistributedGeneralizedNash2022} and signal processing over networks \cite{l.ranDistributedGeneralizedNash2024}. 

\begin{fact}[Basic properties of $\bm{\mathcal{V}}$]
  \label{fact:vGNE_GNEP}
  For the solution sets $\sSGNE{\bm{\mathcal{H}}}{\mathcal{G}}$ in \eqref{eq:GNEP} and $\bm{\mathcal{V}}$ in \eqref{eq:VE}, the following hold:
  \begin{enumerate}[label=\textup{(\roman*)}]
    \item 
    Suppose $\G$ in \eqref{eq:G} is continuous. Then $\bm{\mathcal{V}}$ is closed convex if $\G$ is monotone over $\bm{\mathcal{H}}$ (see, e.g., \cite[Sec. 1.1 and Thm. 2.3.5(a)]{FPvariational2003}). 
    $\bm{\mathcal{V}}$ is nonempty and compact if $\bm{\mathfrak{C}}$ in \eqref{eq:fracture_C} is bounded (see, e.g., \cite[Cor. 2.2.5]{FPvariational2003}). 
    \item (\cite[Thm. 2.1]{facchineiGeneralizedNashGames2007}) 
    $\bm{\mathcal{V}} \subset \sSGNE{\bm{\mathcal{H}}}{\mathcal{G}}$. 
    \item (\cite[Prop. 1.4.2]{FPvariational2003}) 
    For NEP as a simple case of GNEP (see Remark \ref{rem:GNEP_NEP}), 
    $\bm{\mathcal{V}} = \sSGNE{\bm{\mathcal{H}}}{\mathcal{G}}$ holds true. 
    \item 
    In a case where $\f_i\coloneqq\f_0 \ (i\in \mathcal{I})$ with a common differentiable convex function $\f_0:\bm{\mathcal{H}}\to\R$, $\bm{\mathcal{V}} = \textstyle\mathrm{argmin}_{\x \in \bm{\mathfrak{C}}} \ \f_0(\x)$.  
  \end{enumerate}
\end{fact}

In general, the set $\bm{\mathcal{V}}$ in \eqref{eq:VE} is not necessarily singleton. 
This situation induces a challenging equilibrium selection problem: can we design a fair mechanism for each player to reach a certainly desirable v-GNE in $\bm{\mathcal{V}}$ ? 
Related to this question, \cite{g.scutariEquilibriumSelectionPower2012,g.scutariRealComplexMonotone2014,e.benenatiOptimalSelectionTracking2023,w.heDistributedOptimalVariational2024a} proposed a selection model as a hierarchical convex optimization problem (see, e.g., \cite{Yamada-Yamagishi19}), i.e., minimization of a single convex function, called $\uf{}:\bm{\mathcal{H}}\to\R$, over $\bm{\mathcal{V}}$. 
To address this optimization problem, \cite{g.scutariEquilibriumSelectionPower2012,g.scutariRealComplexMonotone2014} 
proposed iterative algorithms of nested structures by introducing an inner loop to solve certain subproblems. 
Quite recently, without using any inner loop, \cite{e.benenatiOptimalSelectionTracking2023,w.heDistributedOptimalVariational2024a} 
proposed to apply the hybrid steepest descent method \cite{Y2001,oguraNonstrictlyConvexMinimization2003,yamadaOHybridQuasi2005,Yamada-Yamagishi19} (see also \cite[Prop. 42]{p.l.combettesFixedPointStrategies2021}) to such~hierarchical~convex~optimization~problems. 

\begin{figure}[htbp]
  \centering
  \includegraphics[clip, width=0.75\columnwidth]{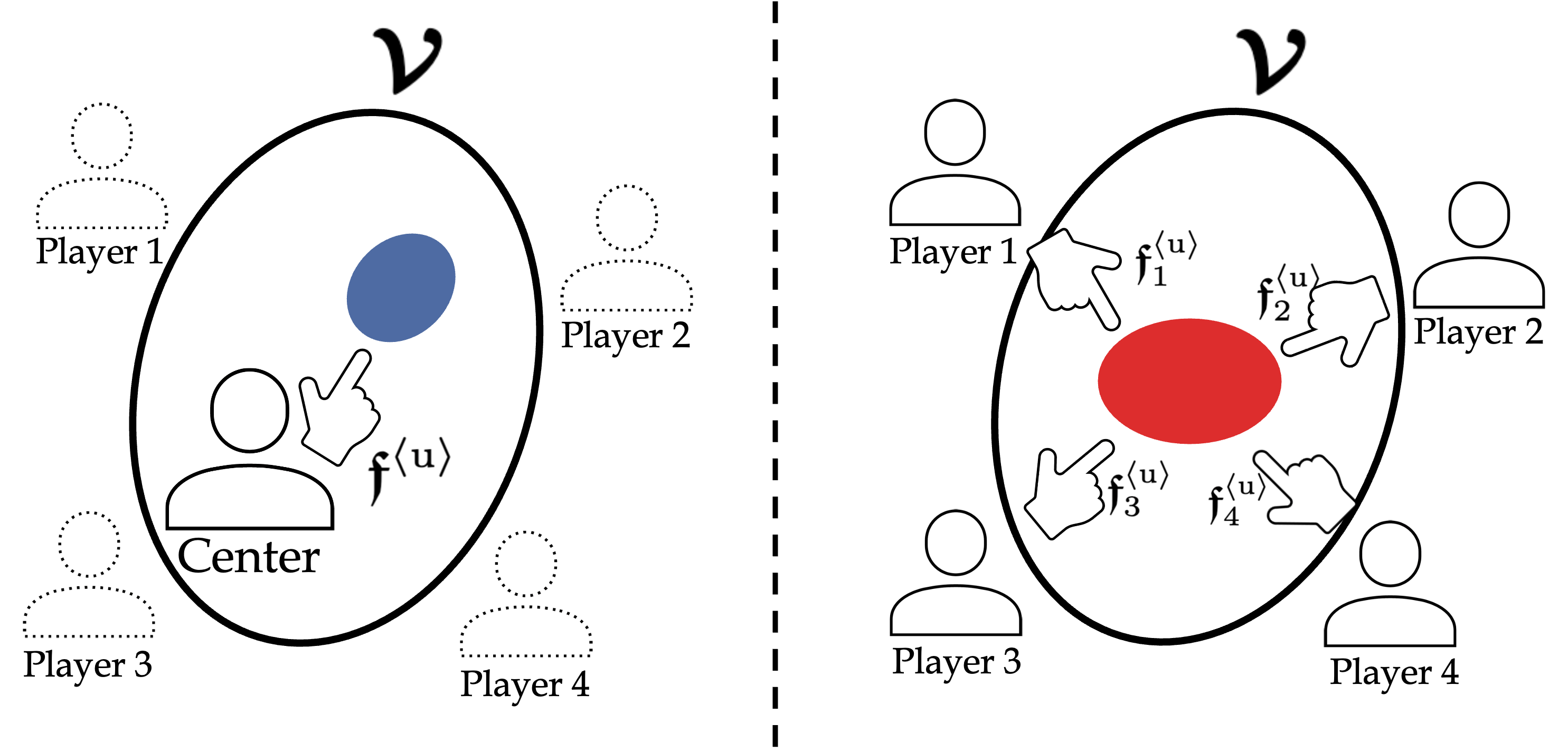}
  \caption{Conceptual comparison of two models for equilibrium selections (existing models $\llbracket$left$\rrbracket$ and proposed model $\llbracket$right$\rrbracket$) over $\bm{\mathcal{V}}$. \\
  \hspace*{1mm} $\llbracket$Left$\rrbracket$ The existing models \cite{g.scutariEquilibriumSelectionPower2012,g.scutariRealComplexMonotone2014,e.benenatiOptimalSelectionTracking2023,w.heDistributedOptimalVariational2024a} have been formulated to choose a special v-GNE by minimizing a single cost function $\uf{}$, designed hopefully by a trusted center, over $\bm{\mathcal{V}}$. \\
  \hspace*{1mm} $\llbracket$Right$\rrbracket$ The proposed model is formulated to choose a special v-GNE, but in a different sense from the existing models, i.e., as an upper-level v-GNE of a new non-cooperative game, over $\bm{\mathcal{V}}$, among all players $i\in\mathcal{I}$ with upper-level cost functions $\uf{i}$ designed by each player $i \in \mathcal{I}$.}
  \label{fig:image_HNEP}
\end{figure}

However, at least from a naive point of view, such hierarchical convex optimization approaches do not seem to achieve ideal fairness among multiple players. 
This is because such approaches certainly require an intervenient, called \emph{center} in this paper, to design a convex function $\uf{}$ (see Fig. \ref{fig:image_HNEP} $\llbracket$Left$\rrbracket$), which deviate from the spirit of non-cooperative game theory pioneered by John Nash. 
In short, another natural question arises: who in the world can design such a function $\uf{}$ certainly according to each player's hope without causing any risk of unexpected bias among players ? 
If we assumed the availability of a center perfectly reliable to all players, we could delegate the authority of designing of $\uf{}$ to the center (the availability of such a trusted center is unfortunately questionable). 

In this paper, by revisiting the spirit of John Nash, we resolve this dilemma without requiring such a trusted center (see Fig. \ref{fig:image_HNEP} $\llbracket$Right$\rrbracket$). 
More precisely, we newly formulate a v-GNEP of hierarchical structure\footnote{
  A similar idea is found in a recent paper \cite{LL2025}, where equilibrium selection is considered for Nash equilibrium problem (see Remark \ref{rem:GNEP_NEP}) as a lower-level noncooperative game. 
} (see Remark \ref{rem:HvGNEP}(i)), where the \emph{lower-level v-GNEP} is given by \eqref{eq:VE} and the \emph{upper-level v-GNEP} is introduced over the solution set $\bm{\mathcal{V}}$ of the lower-level v-GNEP in \eqref{eq:VE}. 
\begin{problem}[Proposed upper-level v-GNEP]
  \label{prob:H_vGNEP}
  Under the setting of the lower-level non-cooperative game formulated in the form of Problem \ref{prob:GNEP}, assume $\bm{\mathcal{V}}$ in \eqref{eq:VE} is a nonempty closed convex set\footnote{\label{fot:sufficient_V} 
  See Fact \ref{fact:vGNE_GNEP}(i) for its sufficient condition.
  }. 
  Then, the upper-level v-GNEP over $\bm{\mathcal{V}}$ is formulated as a variational inequality $\VI(\bm{\mathcal{V}}, \uG)$:
  \begin{equation}
    \label{eq:upper_V}
    \find \x^{\star} \in \uV \coloneqq \bigl\{\x \in \bm{\mathcal{V}} \bigm\vert (\forall \y \in \bm{\mathcal{V}}) \ 
    \product[\bm{\mathcal{H}}]{\uG(\x)}{\y -\x} \geq 0\bigr\}, 
  \end{equation}
  where $\uV \neq \varnothing$ is assumed\footnote{
    The boundedness of $\bm{\mathcal{V}}$ is one of sufficient conditions (see Fact \ref{fact:vGNE_GNEP}(i)).
  },
  \begin{equation}
    \label{eq:uG}
    \bm{\mathfrak{G}} : \bm{\mathcal{H}}\rightarrow\bm{\mathcal{H}} : \x \mapsto \left( \nabla_1 \uf{1}(\x), \dots, \nabla_m \uf{m}(\x) \right), 
  \end{equation}
  for every $i\in\mathcal{I}$, $\uf{i}:\bm{\mathcal{H}}\rightarrow\R$ is the player $i$'s {\em upper-level cost function} such that $\uf{i}(\cdot;\x_{\smallsetminus i}):\mathcal{H}_i \rightarrow \mathbb{R}$ is convex and differentiable for every $\x \in \bm{\mathcal{H}}$, and $\nabla_i \uf{i}(\cdot;\x_{\smallsetminus i})$ is the gradient of $\uf{i}(\cdot;\x_{\smallsetminus i})$ (see Remark \ref{rem:HvGNEP}(iii) regarding the design of $\uf{i}$).
\end{problem}
\begin{remark}[On Problem \ref{prob:H_vGNEP}]\label{rem:HvGNEP}\ 
  \begin{enumerate}[label=\textup{(\roman*)}]
    \item (Interpretation of Problem \ref{prob:H_vGNEP}). 
    For the solution set $\bm{\mathcal{V}}$ of the lower-level v-GNEP $\VI(\bm{\mathfrak{C}}, \G)$, 
    $\SGNE{\bm{\mathcal{H}}}{\mathcal{G}}{\uf{i}}{\mathcal{H}_i}{\Id}{\bm{\mathcal{V}}}$ is well-defined, in a similar way to Problem \ref{prob:GNEP}, as 
    \begin{align}
      \label{eq:HGNEP_V}
        &\find \x \in \SGNE{\bm{\mathcal{H}}}{\mathcal{G}}{\uf{i}}{\mathcal{H}_i}{\Id}{\bm{\mathcal{V}}} = \Bigl\{(\bar{x}_1, \dots, \bar{x}_m) \in \bm{\mathcal{H}} \bigm\vert \\
        &\qquad (\forall i \in \mathcal{I}) \ \bar{x}_i \in \argmin[x_i \in \mathcal{H}_i] \ \uf{i}(x_i; \bar{\x}_{\smallsetminus i}) \ \st (x_i; \bar{\x}_{\smallsetminus i}) \in \bm{\mathcal{V}} \Bigr\} \subset \bm{\mathcal{V}}. 
    \end{align}
    The upper-level v-GNEP in \eqref{eq:upper_V} is the v-GNEP for $\SGNE{\bm{\mathcal{H}}}{\mathcal{G}}{\uf{i}}{\mathcal{H}_i}{\Id}{\bm{\mathcal{V}}}$. 
    By Fact \ref{fact:vGNE_GNEP}(ii), we have
    \begin{equation}
      \label{eq:V_S_V_S}
      \uV \subset \SGNE{\bm{\mathcal{H}}}{\mathcal{G}}{\uf{i}}{\mathcal{H}_i}{\Id}{\bm{\mathcal{V}}} 
      \subset \bm{\mathcal{V}}
      \subset \SGNE{\bm{\mathcal{H}}}{\mathcal{G}}{\f_i}{C_i}{\Li}{D}. 
    \end{equation}
    \item (Challenge in Problem \ref{prob:H_vGNEP}). 
    In general, $\bm{\mathcal{V}}$ in Problems \eqref{eq:upper_V} and \eqref{eq:HGNEP_V} is possibly infinite set. 
    So far, even for just finding a point in $\bm{\mathcal{V}}$ without any aimed selection, a certain iterative approximation, via a strategically produced infinite sequence, has been required (see, e.g., \cite{buiNashEquilibrium2021,g.belgioiosoDistributedGeneralizedNash2022}). 
    This situation suggests that, \emph{in order to achieve the mission of Problem \ref{prob:H_vGNEP}, we have to elaborate a mathematical strategy for meeting, as much as possible, all players' upper-level requirements $\uf{i} \ (i\in\mathcal{I})$ without stopping the analysis of the whole view of $\bm{\mathcal{V}}$}. 
    \item 
    $\uf{i}:\bm{\mathcal{H}}\rightarrow\R$ can be designed, by each player $i(\in\mathcal{I})$, independently of the lower-level non-cooperative game.
    \item By assigning a common differentiable convex function $\uf{0}:\bm{\mathcal{H}}\to\R$ to $\uf{i} \ (\forall i \in \mathcal{I})$, the upper-level v-GNEP in \eqref{eq:upper_V} reproduces the hierarchical convex optimization problems over $\bm{\mathcal{V}}$ \cite{g.scutariEquilibriumSelectionPower2012,g.scutariRealComplexMonotone2014,e.benenatiOptimalSelectionTracking2023,w.heDistributedOptimalVariational2024a}:
    \begin{equation}
      \label{eq:upper_V_common}
      \find \x^{\star} \in  \argmin[\x \in \bm{\mathcal{V}}] \ \uf{0}(\x).
    \end{equation}
  \end{enumerate}
\end{remark}

For finding a solution of Problem \ref{prob:H_vGNEP}, we propose an iterative algorithm (see Theorem~\ref{thm:modified_HSDM}). 
Inspired by \cite{briceno-ariasMonotoneOperatorMethods2013,e.benenatiOptimalSelectionTracking2023}, the proposed algorithm is designed as an application of the hybrid steepest descent method \cite{yamadaOHybridQuasi2005} to a fixed point expression of $\bm{\mathcal{V}}$ via the so-called forward-backward-forward operator \cite{tsengModifiedForwardBackwardSplitting2000} which is quasi-nonexpansive (see Proposition \ref{prop:FBF_fix}(ii)). 

To show a clear distinction of the proposed equilibrium selection achievable by Problem \ref{prob:H_vGNEP} from the existing equilibrium selection achievable by \eqref{eq:upper_V_common}, we present Problem \ref{prob:cycle_selection} as an instance of Problem \ref{prob:H_vGNEP}. 
Problem \ref{prob:cycle_selection} is formulated as a v-GNEP for finding a \emph{cycle} \cite{BAILLON2012400,p.l.combettesFixedPointStrategies2021,OJMO_2021__2__A5_0} over $\bm{\mathcal{V}}$ and includes at least~one~example~that~can~not~be characterized~as~an~instance~of~the~existing~problem~\eqref{eq:upper_V_common}. 

The remainder of this paper is organized as follows. 
Section \ref{sec:preliminaries} introduces (i) selected tools in fixed point theory of quasi-nonexpansive operator, and (ii) cycles in the cyclic projection algorithm. 
In Section \ref{sec:Algorithm}, we present the proposed algorithm for Problem \ref{prob:H_vGNEP} and an example of the equilibrium selection based on cycles. 
Section \ref{sec:experiments} presents numerical experiments, followed by conclusion in Section \ref{sec:conclusion}. 
A preliminary short version of this paper was presented at a conference \cite{Matsuo2025}. 

\noindent\textbf{Notation} \ 
$\mathbb{N}$ and $\mathbb{R}$ denote respectively the set of all nonnegative integers and the set of all real numbers.
Let$(\mathcal{H}, \product[\mathcal{H}]{\cdot}{\cdot}, \|\cdot\|_{\mathcal{H}})$ and $(\mathcal{K}, \product[\mathcal{K}]{\cdot}{\cdot}, \|\cdot\|_{\mathcal{K}})$ be finite dimensional real Hilbert spaces.
The set of all bounded linear operators from $\mathcal{H}$ to $\mathcal{K}$ is denoted by $\mathcal{B}(\mathcal{H},\mathcal{K})$.
$\Id\in\mathcal{B}(\mathcal{H},\mathcal{H})$ and $O_{\mathcal{H}}\in\mathcal{B}(\mathcal{H},\mathcal{H})$ denote respectively the identity operator and the zero operator. 
For $\mathcal{L}\in\mathcal{B}(\mathcal{H},\mathcal{K})$, $\mathcal{L}^* \in\mathcal{B}(\mathcal{K},\mathcal{H})$ denotes the adjoint operator of $\mathcal{L}$, i.e., $(\forall (x,y)\in\mathcal{H}\times\mathcal{K})$ $\langle \mathcal{L}x, y\rangle_{\mathcal{K}} = \langle x, \mathcal{L}^*y\rangle_{\mathcal{H}}$. 
The range and the kernel of $\mathcal{L}\in\mathcal{B}(\mathcal{H},\mathcal{K})$ are denoted respectively by $\mathrm{ran}(\mathcal{L})=\{\mathcal{L}x\in\mathcal{K}\mid x\in\mathcal{H}\}$ and by $\operatorname{ker}(\mathcal{L}) = \{x \in \mathcal{H} \mid \mathcal{L}x = 0_{\mathcal{K}}\}$. 
The operator norm of $\mathcal{L}\in\mathcal{B}(\mathcal{H},\mathcal{K})$ is denoted by $\|\mathcal{L}\|_{\mathrm{op}}\coloneqq\sup_{x\in\mathcal{H},\|x\|_{\mathcal{H}}\leq1}\|\mathcal{L}x\|_{\mathcal{K}}$.
For sets $S_1\subset\mathcal{H}$ and $S_2 \subset \mathcal{K}$, the image of $S_1$ and the preimage of $S_2$ under $\mathcal{L}\in\mathcal{B}(\mathcal{H},\mathcal{K})$ are denoted respectively by $\mathcal{L}(S_1):=\{\mathcal{L}x \in \mathcal{K} \vert x \in S_1\}$ and by $\mathcal{L}^{-1}(S_2):=\{x \in \mathcal{H} \vert \mathcal{L}x \in S_2\}$. 
For a set $S\subset\mathcal{H}$,
(i) the linear span of $S$ is
$\operatorname{span}(S)\coloneqq \bigl\{\sum_{i=1}^{n}\alpha_{i} x_{i}\in\mathcal{H}\mid n\in\mathbb{N}, \alpha_{i}\in\mathbb{R},x_{i}\in S\bigr\}$,
(ii) the affine hull of $S$ is
$\operatorname{aff}(S)\coloneqq x+\operatorname{span}(S-x)$, where $x\in S$, $S-x\coloneqq \{w-x\in\mathcal{H}\mid w\in S\}$, and $x+\operatorname{span}(S-x)=y+\operatorname{span}(S-y)$ holds true for all $x,y\in S$, and 
(iii) the interior of $S$ is
$\operatorname{int}(S)\coloneqq \{x\in S\mid (\exists\rho>0)\,B(x;\rho)\subset S\}$,
where $B(x;\rho)\coloneqq\{y\in\mathcal{H}\mid\|y-x\|_{\mathcal{H}}<\rho\}$ is an open ball with center $x\in\mathcal{H}$ and radius $\rho>0$.
The relative interior of a convex set $S\subset\mathcal{H}$ is $\operatorname{ri}(S)\coloneqq\{x\in S\mid (\exists\rho>0)\,B(x;\rho)\cap \operatorname{aff}(S)\subset S\}$.
The power set of $\mathcal{H}$, denoted by $2^{\mathcal{H}}$, is the collection of all subsets of $\mathcal{H}$, i.e., $2^{\mathcal{H}}\coloneqq\{S\mid S\subset\mathcal{H}\}$. 

A set $S\subset\mathcal{H}$ is said to be convex if $(1-\theta)x+\theta y\in S$ for all $(x,y,\theta)\in S\times S\times [0,1]$. 
A function $f:\mathcal{H}\to(-\infty, +\infty]$ is said to be proper if $\dom(f):= \{x \in \mathcal{H} \mid f(x) < +\infty\}\neq\varnothing$, lower semicontinuous if $\mathrm{lev}_{\leq \alpha}(f):=\{x\in\mathcal{H}\mid f(x) \leq \alpha\}$ is closed for every $\alpha \in \R$, and convex if $f(\theta x + (1-\theta)y) \leq \theta f(x) + (1-\theta)f(y)$ for all $(x,y,\theta)\in \operatorname{dom}(f)\times \operatorname{dom}(f) \times (0,1)$, respectively. 
The set of all proper lower semicontinuous convex functions defined over $\mathcal{H}$ is denoted by $\Gamma_0(\mathcal{H})$. 
Let $f\in\Gamma_{0}(\mathcal{H})$. 
The conjugate (also named Legendre-Fenchel transform) of $f$ is defined by $f^*:\mathcal{H}\to (-\infty, \infty]:u\mapsto\sup_{x\in\mathcal{H}}[\langle x, u\rangle_{\mathcal{H}} - f(x)]$. 
The subdifferential of $f$ is defined as the set-valued operator $\partial f:\mathcal{H}\to 2^{\mathcal{H}}:x\mapsto\{u\in\mathcal{H}\mid \langle y-x,u \rangle_{\mathcal{H}}+f(x)\leq f(y),\forall y\in\mathcal{H}\}$. 
Let $C\subset\mathcal{H}$ be a nonempty closed convex set. 
The indicator function of $C$ is defined as $\iota_C(x)=0$ if $x \in C$ and $\iota_C(x)=+\infty$ if $x \notin C$ (Note: $\iota_C \in \Gamma_0(\mathcal{H})$). 
The metric projection onto $C$ is defined by $P_C:\mathcal{H}\to\mathcal{H}:x\mapsto\textstyle\mathrm{argmin}_{y\in C}\|x - y\|_{\mathcal{H}}$. 
The distance to $C$ is defined by $d(\cdot, C): \mathcal{H} \ni x \mapsto \min_{y \in C} \|x - y\|_{\mathcal{H}}$. 
An operator $\mathcal{A}:\mathcal{H}\to\mathcal{H}$ is said to be (i) monotone over $S\subset \mathcal{H}$ if $\product[\mathcal{H}]{\mathcal{A}(x) - \mathcal{A}(y)}{x - y} \geq 0$ for all $x, y \in S$, paramonotone over $S$ if $\mathcal{A}$ is monotone over $S$ and for all $x, y \in S$, $\product[\mathcal{H}]{\mathcal{A}(x) - \mathcal{A}(y)}{x - y} = 0 \Leftrightarrow \mathcal{A}(x) = \mathcal{A}(y)$, and (ii) Lipschitz continuous with a Lipschitz constant $\kappa>0$ (or $\kappa$-Lipschitzian) over $S$ if $\|\mathcal{A}(x) - \mathcal{A}(y)\|_{\mathcal{H}} \leq \kappa\|x - y\|_{\mathcal{H}}$ for all $x, y \in S$. 

%% file: docs/Preliminaries.tex
\section{Preliminaries}
\label{sec:preliminaries}
\subsection{Quasi-Nonexpansive Operator and Hybrid Steepest Descent Method}
\begin{definition}[{Quasi-nonexpansive operators, e.g., \cite[Sec. B]{yamadaOHybridQuasi2005}}, \cite{yamadaMinimizingMoreauEnvelope2011}]
  \label{df:quasi_nonexpansive} \ 
  \begin{enumerate}[label=\textnormal{(\roman*)}]
    \item (Nonexpansive operator). An operator $T:\mathcal{H}\to\mathcal{H}$ is said to be \emph{nonexpansive} if $T$ is $1$-Lipschitzian over $\mathcal{H}$, i.e., 
    \begin{equation}
      (\forall x, y \in \mathcal{H}) \ \|T(x) - T(y)\|_{\mathcal{H}} \leq \|x - y\|_{\mathcal{H}}.
    \end{equation}
    For example, $P_C$ is nonexpansive for a nonempty closed convex set $C\subset\mathcal{H}$. 
    \item (Quasi-nonexpansive operator). 
    An operator $T$ satisfying $\Fix(T):=\{x \in \mathcal{H} \mid T(x) = x\} \neq \varnothing$ is said to be \emph{quasi-nonexpansive} if 
    \begin{equation}
      \bigl(\forall x \in \mathcal{H}, \forall z \in \Fix(T)\bigr) \ \|T(x) - z\|_{\mathcal{H}} \leq \|x - z\|_{\mathcal{H}}. 
    \end{equation}
    For example, a nonexpansive operator $T:\mathcal{H}\to\mathcal{H}$ with $\Fix(T)\neq\varnothing$ is quasi-nonexpansive. 
    \item (Attracting operator). 
    A quasi-nonexpansive operator $T:\mathcal{H}\to\mathcal{H}$ is said to be \emph{attracting} if 
    \begin{equation}
      \bigl(\forall x \notin \Fix(T), \forall z \in \Fix(T)\bigr) \ \|T(x) - z\|_{\mathcal{H}} < \|x - z\|_{\mathcal{H}}. 
    \end{equation}
    \item (Strongly attracting operator). 
    An attracting operator $T$ is said to be ($\eta$-)\emph{strongly} attracting\footnote{Such $T$ is also said to be \emph{strongly quasi-nonexpansive} in \cite{cegielskiAlgorithmSolvingVariational2013}.} if 
    \begin{equation}
      (\exists \eta> 0, \forall x \in \mathcal{H}, \forall z \in \Fix(T)) \ \|T(x) - z\|^2_{\mathcal{H}} \leq \|x - z\|^2_{\mathcal{H}} - \eta\|T(x) - x\|^2_{\mathcal{H}}. 
    \end{equation}
    \item (Averaged operator). 
    A quasi-nonexpansive operator $T:\mathcal{H}\to\mathcal{H}$ is said to be ($\alpha$-)\emph{averaged} if there exist some $\alpha \in (0, 1)$ and some quasi-nonexpansive operator $U : \mathcal{H} \rightarrow \mathcal{H}$ such that 
    \begin{equation}
      T = (1 - \alpha)\Id + \alpha U. 
    \end{equation}
    In this case, $\Fix(T) = \Fix(U) \neq \varnothing$ holds true. 
  \end{enumerate}
\end{definition}

\begin{fact}[{Selected properties of quasi-nonexpansive operators, e.g., \cite[Prop. 1]{yamadaOHybridQuasi2005}}]\label{fact:quasi_nonexpansive} \ 
  \begin{enumerate}
    \renewcommand{\labelenumi}{\textnormal{(\roman{enumi})}}
    \item For a quasi-nonexpansive operator $T:\mathcal{H}\rightarrow\mathcal{H}$, $\Fix(T)$ is closed convex. 
    \item For $\alpha \in (0, 1)$ and a quasi-nonexpansive operator $T:\mathcal{H}\rightarrow\mathcal{H}$, $T$ is $\alpha$-averaged if and only if $T$ is $\frac{1-\alpha}{\alpha}$-strongly attracting. 
    \item Let $T_1: \mathcal{H}\rightarrow \mathcal{H}$ be quasi-nonexpansive and 
    $T_2: \mathcal{H}\rightarrow \mathcal{H}$ be attracting quasi-nonexpansive with 
    $\Fix(T_1) \cap \Fix(T_2) \neq \varnothing$.
    Then $T_2 \circ T_1:\mathcal{H}\rightarrow\mathcal{H}$ is a quasi-nonexpansive mapping with $\Fix(T_2 \circ T_1) = \Fix(T_1) \cap \Fix(T_2)$. 
    Moreover, if $T_1$ is $\alpha_1$-strongly attracting and $T_2$ is $\alpha_2$-strongly attracting, then $T_2 \circ T_1$ is $\frac{\alpha_1 \alpha_2}{\alpha_1 + \alpha_2}$-strongly attracting. 
  \end{enumerate} 
\end{fact}

\begin{definition}[{Quasi-shrinking operator \cite[Def. 1]{yamadaOHybridQuasi2005}}]
  \label{df:quasi_shrinking}
  Let $T:\mathcal{H}\rightarrow\mathcal{H}$ be a quasi-nonexpansive operator satisfying $\Fix(T) \cap C \neq \varnothing$ for some closed convex set $C \subset \mathcal{H}$. 
  The operator $T$ is said to be \emph{quasi-shrinking} on $C$ if 
  \begin{align}
    D: [0, \infty) \ni r \mapsto 
    \begin{cases}
      \displaystyle \inf_{x \in C \cap \bigl(\operatorname{lev}_{\geq r} (d_{\Fix(T)})\bigr)} d\bigl(x, \Fix(T)\bigr) - d\bigl(T(x), \Fix(T)\bigr) \\
      \qquad \ \mathrm{if} \ C \cap \bigl(\operatorname{lev}_{\geq r} (d_{\Fix(T)})\bigr) \neq \varnothing, \\
      \infty \quad \ \mathrm{otherwise}
    \end{cases}
  \end{align}
  satisfies $D(r) = 0 \Leftrightarrow r = 0$, where $d_{\Fix(T)}:\mathcal{H} \to [0,\infty): x\mapsto d(x, \Fix(T))$ and $\operatorname{lev}_{\geq r} (d_{\Fix(T)})\coloneqq \{x \in \mathcal{H} \vert d_{\Fix(T)}(x) \geq r\}$ (Note: ($\Leftarrow$) always holds true). 
\end{definition}

\begin{fact}[{A sufficient condition 
  to be quasi-shrinking \cite[Prop. 2.11]{cegielskiAlgorithmSolvingVariational2013}}]
  \label{fact:quasi_shrinking_stquasi}
  Let $C \subset \mathcal{H}$ be a bounded closed convex set and let $T:\mathcal{H} \rightarrow \mathcal{H}$ a quasi-nonexpansive operator with $\Fix(T) \cap C \neq \varnothing$. 
  Suppose that 
  (i) $T:\mathcal{H} \rightarrow \mathcal{H}$ is strongly attracting, and
  (ii) $T - \Id$ is demi-closed at $\bm{0}_{\mathcal{H}}$, i.e., for any sequence $(x_n)_{n\in\N} \subset \mathcal{H}$, 
  \begin{equation}
    \bigl((\exists x\in\mathcal{H})\lim_{n\to\infty}x_n = x  \ \mathrm{and} \lim_{n\to\infty}\|T(x_n) - x_n\| = 0\bigr) \Rightarrow \underbrace{T(x) - x = \bm{0}_{\mathcal{H}}}_{\Leftrightarrow x \in \Fix(T)}. 
  \end{equation}
  Then $T$ is quasi-shrinking on $C$. 
\end{fact}

Consider the variational inequality $\VI(\Fix(T), F)$:
\begin{equation}
  \label{eq:VI_fix}
  \find x \in \Fix(T) \ \st \bigl(\forall y \in \Fix(T)\bigr) \ \product[]{F(x)}{y-x} \geq 0,
\end{equation}
where $T:\mathcal{H}\to\mathcal{H}$ is a quasi-nonexpansive operator with $\Fix(T) \neq \varnothing$ and $F:\mathcal{H}\to\mathcal{H}$ is a monotone operator over $\Fix(T)$. 
The sequence $(x_n)_{n \in \N}$ generated by the hybrid steepest descent method \cite{yamadaOHybridQuasi2005} 
\begin{equation}
  \label{eq:HSDM_preliminary}
  (n \in \N) \ x_{n+1} = T(x_n) - \lambda_{n+1} F\bigl( T(x_n) \bigr)
\end{equation}
successively approximates a solution of the problem \eqref{eq:VI_fix} 
under the conditions in the next fact. 
\begin{fact}[{Hybrid steepest descent method for quasi-nonexpansive operators \cite[Theorem 5]{yamadaOHybridQuasi2005}}]
  \label{fact:HSDM}
  Suppose that 
  (a) $T:\mathcal{H} \rightarrow \mathcal{H}$ is quasi-nonexpansive with bounded $\Fix(T) \neq \varnothing$, 
  (b) $F:\mathcal{H} \rightarrow \mathcal{H}$ is paramonotone over $\Fix(T)$ and Lipschitz continuous over $T(\mathcal{H}) \coloneqq \{T(x) \in \mathcal{H} \vert x \in \mathcal{H}\}$. 
  Set~\footnote{By Fact \ref{fact:quasi_nonexpansive}(i), $\Fix(T) \neq \varnothing$ is closed convex. Since (i) $\Fix(T)$ is bounded and (ii) $F$ is paramonotone over $\Fix(T)$ and Lipschitzian over $T(\mathcal{H})$ by the assumption, $\Gamma$ is nonempty closed convex by \cite[Sec. 1.1, Cor. 2.2.5, Thm. 2.3.5(a)]{FPvariational2003}.} $\Gamma \coloneqq \{x \in \Fix(T)\mid \product[]{F(x)}{y - x} \geq 0, \forall y \in \Fix(T)\}$. 
  By using $(\lambda_n)_{n\in\N} \subset [0, +\infty)$ such that (H1) $\lim_{n\rightarrow\infty}\lambda_n = 0$ and (H2) $\sum_{n=1}^{\infty}\lambda_n = \infty$, for any $x_0 \in \mathcal{H}$, the sequence $(x_n)_{n \in \N}$ generated by \eqref{eq:HSDM_preliminary} satisfies $\lim_{n\rightarrow\infty}d(x_n, \Gamma) = 0$ 
  if there exists nonempty bounded closed convex set $C\subset\mathcal{H}$ satisfying $(x_n)_{n\in\N} \subset C$ and $T$ is quasi-shrinking on $C$. 
\end{fact}

\begin{remark}[On the hybrid steepest descent method]
  \label{rem:HSDM}
  Under the setting of Fact \ref{fact:HSDM}, we have:
  \begin{enumerate}
    \renewcommand{\labelenumi}{\textnormal{(\roman{enumi})}}
    \item There exists\footnote{Since $(x_n)_{n\in\N}$ is bounded by the assumption, there exists a subsequence $(x_{n_k})_{k\in\N}$ that converges to some $\tilde{x} \in \mathcal{H}$. By $\lim_{n\rightarrow\infty}d(x_n, \Gamma) = 0$ and the continuity of $d(\cdot, \Gamma)$, we have $\lim_{k\rightarrow\infty}d(x_{n_k}, \Gamma) = d(\tilde{x}, \Gamma) = 0$, which implies $\tilde{x}\in\Gamma$. } a cluster point of $(x_n)_{n\in\N}$ generated by \eqref{eq:HSDM_preliminary}, and any cluster point of $(x_n)_{n\in\N}$ belongs to $\Gamma$. 
    \item Let $\Phi:\mathcal{H}\to\R$ be a convex and differentiable function. 
    In a case of the problem \eqref{eq:VI_fix} where $F$ is chosen as $F = \nabla\Phi$, \eqref{eq:VI_fix} reproduces a hierarchical convex optimization problem over $\Fix(T)$:
    \begin{equation}
      \label{eq:hierarchical_convex}
      \find \bar{x} \in \argmin[x\in\Fix(T)] \Phi(x).
    \end{equation}
    For applications of the hybrid steepest descent method to Problem \eqref{eq:hierarchical_convex}, see, e.g., \cite{yamadaMinimizingMoreauEnvelope2011,Yamada-Yamagishi19}. 
  \end{enumerate}
\end{remark}

\subsection{Cycles as Nash Equilibria not Characterizable via Optimization Problem}
\label{subsec:preriminary_cycle}
We introduce the notion of cycles, which has served as a useful analytic tool for inconsistent convex feasibility problems (see, e.g., \cite{censor2018algorithms,p.l.combettesFixedPointStrategies2021,OJMO_2021__2__A5_0} and reference therein). 
\begin{definition}[Cycles in the cyclic projection algorithm]\label{def:cycle}
    Let $m$ be an integer at least equal to $2$ and let $(K_1, \dots, K_m)$ be an ordered family of nonempty closed convex subsets of $\mathcal{H}$. 
    Then, a tuple $(\bar{x}_1, \dots, \bar{x}_m) \in \mathcal{H}^m$ is said to be a \emph{cycle} associated with $(K_1, \dots, K_m)$ if 
    \begin{equation}
      \label{eq:cycle_POCS}
      \begin{aligned}
        \bar{x}_2 = P_{K_2}(\bar{x}_1), \bar{x}_3 = P_{K_3}(\bar{x}_2), \dots, \bar{x}_{m} = P_{K_{m}}(\bar{x}_{m-1}), \bar{x}_1 = P_{K_1}(\bar{x}_m). 
      \end{aligned}
    \end{equation}
    The set of all cycles associated with $(K_1, \dots, K_m)$ is denoted by $\mathrm{cyc}(K_1, \dots, K_m)$. 
\end{definition}
\begin{remark}[POCS algorithm]
  \label{rem:POCS}
  If at least one of $\{K_1, \dots, K_m\}$ is bounded, then $\mathrm{cyc}(K_1, \dots, K_m) \neq \varnothing$ is guaranteed \cite{Gubin1967TheMO}. 
  Moreover, if $P_{K_i}$ is available as a computable operator for every $i \in \{1, \dots, m\}$, we can apply the \emph{Projection Onto Convex Sets (POCS) algorithm} \cite{BLEGMAN1965,EREMIN1965} for finding a cycle \eqref{eq:cycle_POCS}. 
\end{remark}

For $m=2$ case in \eqref{eq:cycle_POCS}, the set $\mathrm{cyc}(K_1, K_2)$ of all cycles can be expressed as the solution set of the following optimization problem \cite{cheneyProximityMapsConvex1959,bauschkeConvergenceNeumannsAlternating1993,BAUSCHKE1994418}:
\begin{equation}
  \label{eq:best_approximation}
  \mathrm{cyc}(K_1, K_2) = \argmin[x_1 \in K_1, x_2 \in K_2] \|x_1 - x_2\|. 
\end{equation}

However, for $m\geq 3$ cases, the set of all cycles can not be characterized as optimization problems. 
\begin{fact}[{There is no variational characterization of the cycles \cite[Thm. 2.3]{BAILLON2012400}}]
  \label{fact:no_variational_cycle}
  Suppose that $\dim(\mathcal{H}) \geq 2$ and let $m$ be an integer at least equal to $3$. 
  There exists no function $\Phi : \mathcal{H}^m\to\R $ such that, for every ordered family of nonempty closed convex subsets $(K_1, \dots, K_m)$ of $\mathcal{H}$, 
  \begin{equation}
    \label{eq:no_variational_cycle}
    \mathrm{cyc}(K_1, \dots, K_m) = \argmin[x_1\in K_1, \dots, x_m \in K_m] \Phi(x_1, \dots, x_m). 
  \end{equation}
  In other words, we can not construct any function $\Phi:\mathcal{H}^m\to\R$ satisfying \eqref{eq:no_variational_cycle} in a unified way, i.e., independently of the choices of ordered family of nonempty closed convex sets $(K_1, \dots, K_m)$. 
\end{fact}

Meanwhile, as discussed in \cite[Exm. 9.1.4]{briceno-ariasMonotoneOperatorMethods2013} and \cite[Sec. VI.A]{p.l.combettesFixedPointStrategies2021}, for every ordered family of nonempty closed convex sets $(K_1, \dots, K_m)$, $\mathrm{cyc}(K_1, \dots, K_m)$ can be expressed\footnote{\label{fot:cycle_game_eq}We can check that $\mathrm{cyc}(K_1, \dots, K_m)$ is the solution set of \eqref{eq:cycle_game} as follows. 
\begin{equation}
  \label{eq:cycle_game_confirm}
  \begin{aligned}
    \mathrm{cyc}(K_1, \dots, K_m) 
    &= \{(\bar{x}_1, \dots, \bar{x}_m) \in \mathcal{H}^m \vert (\forall i \in \mathcal{I}) \ \bar{x}_i = P_{K_i}(\bar{x}_{i-1})\} \\
    &= \{(\bar{x}_1, \dots, \bar{x}_m) \in \mathcal{H}^m \vert (\forall i \in \mathcal{I}) \ \bar{x}_i \in \textstyle \mathrm{argmin}_{x_i \in K_i}  \frac{1}{2}\|x_i - \bar{x}_{i-1}\|^2\}. 
  \end{aligned}
\end{equation}
}, in a unified way, as the solution set 
\begin{equation}
  \mathrm{cyc}(K_1, \dots, K_m) = \SGNE{\mathcal{H}^m}{\mathcal{H}^m}{\f_i}{K_i}{\Id}{\mathcal{H}^m}, \f_i(\x)\coloneqq\textstyle\frac{1}{2}\|x_i - x_{i-1}\|^2
\end{equation}
of the following Nash equilibrium problem with the convention $0 = m$:
\begin{align}
  \label{eq:cycle_game}
  \find (\bar{x}_1, \dots, \bar{x}_m) \in \mathcal{H}^m \ \st (\forall i \in \mathcal{I}) \ 
  \bar{x}_i \in \argmin[x_i \in K_i] \ \frac{1}{2}\|x_i - \bar{x}_{i-1}\|^2. 
\end{align}

%% file: docs/Algorithm.tex
\section{Proposed Algorithm and Application}
\label{sec:Algorithm}
\subsection{Reformulation of Upper-level v-GNEP as a Variational Inequality over Fixed Point Set of Quasi-Nonexpansive Operator}
Fixed point theory has been offering powerful ideas for solving the GNEP \cite{p.l.combettesFixedPointStrategies2021,g.belgioiosoDistributedGeneralizedNash2022}. 
We propose to use such ideas for solving Problem \ref{prob:H_vGNEP} under Assumption \ref{asm:G_monotone_Lipschitzian} below. 
\begin{assumption}
  \label{asm:G_monotone_Lipschitzian}
  Under the setting of Problem \ref{prob:H_vGNEP}, assume that:
  \begin{enumerate}[label=\textnormal{(\alph*)}]
    \item \label{item:asm_sum_chain} For $\Li:\bm{\mathcal{H}}\to\mathcal{G}, \C \subset \bm{\mathcal{H}}$, and $D \subset \mathcal{G}$ in Problem \ref{prob:GNEP}\textit{\ref{item:Li_C_D}}, $\partial \bigl( \iota_{\bm{C}} + \iota_D \circ \Li \bigr) = \partial \iota_{\bm{C}} + \Li^{\ast} \circ \partial \iota_D \circ \Li \footnote{For its sufficient conditions, see Remark \ref{rem:bound_cond} and \cite[Thm. 16.47]{BC2017}.}$.
    \item \label{item:asm_G} $\G:\bm{\mathcal{H}}\to\bm{\mathcal{H}}$ in \eqref{eq:G} is $\kappa_{\G}(>0)$-Lipschitzian and monotone over $\bm{\mathcal{H}}$ [this assumption ensures closed convexity of $\bm{\mathcal{V}}$ in \eqref{eq:VE} (see Fact \ref{fact:vGNE_GNEP}(i))]. 
    \item \label{item:asm_frakG} $\bm{\mathfrak{G}}:\bm{\mathcal{H}}\rightarrow\bm{\mathcal{H}}$ in \eqref{eq:uG} is $\kappa_{\bm{\mathfrak{G}}}(>0)$-Lipschitzian and paramonotone over $\bm{\mathcal{H}}$~\footnote{For its sufficient conditions, see, e.g., \cite[Chap. 22]{BC2017}. In particular, for a simple case where $\uf{i} \coloneqq \uf{0} \ (i\in\mathcal{I})$ with a common differentiable convex function $\uf{0}:\bm{\mathcal{H}}\to\R$ (see also Remark \ref{rem:HvGNEP}(iv)), $\uG$ can be expressed as $\uG = \nabla \uf{0}$, and thus the paramonotonicity of $\uG$ is automatically guaranteed \cite[Exm. 22.4]{BC2017}. } [this assumption ensures closed convexity of $\uV$ in \eqref{eq:upper_V} (see Fact \ref{fact:vGNE_GNEP}(i))].  
  \end{enumerate}
\end{assumption}
The next proposition is motivated by \cite{Franci2020FBFforGeneralizedNash,e.benenatiOptimalSelectionTracking2023}. 
\begin{proposition}[Fixed point expression of v-GNE via Forward-\linebreak Backward-Forward operator]
  \label{prop:FBF_fix}
  Under Assumption \ref{asm:G_monotone_Lipschitzian}, define $\TFBF:\bm{\mathcal{H}} \times\mathcal{G} \rightarrow \bm{\mathcal{H}} \times\mathcal{G}$ \footnote{
    Since $\B$ is maximally monotone (see Definition \ref{def:monotone_operator} and Lemma \ref{lemma:A_B_monotone}), $(\Id + \gamma \B)^{-1}$ is single-valued (see Fact \ref{fact:maximally_monotone_operator}\ref{item:single_resolvent}). Hence, $\T_{\mathrm{FB}}$ and $\TFBF$ in \eqref{eq:FBF} are also single-valued.
  }, 
  with $\gamma \in \bigl(0, 1/(\kappa_{\G} + \|\Li\|_{\mathrm{op}})\bigr)$, by 
  \begin{equation}
    \label{eq:FBF}
    \TFBF := (\Id - \gamma\A) \circ \underbrace{(\Id + \gamma \B)^{-1} \circ (\Id - \gamma \A)}_{=:\T_{\FB}} + \gamma \A, 
  \end{equation}
  and its $\alpha(\in(0, 1))$-averaged version by 
  \begin{equation}
    \label{eq:Talpha}
    \T_{\FBF}^{\alpha} := (1-\alpha)\Id + \alpha \TFBF, 
  \end{equation}
  where
  \begin{align}
    \label{eq:A}
    &\A:\bm{\mathcal{H}} \times\mathcal{G} \rightarrow \bm{\mathcal{H}} \times\mathcal{G} : (\x, u) \mapsto (\G(\x)+\Li^*u, -\Li\x), \\
    \label{eq:B}
    &\B:\bm{\mathcal{H}} \times\mathcal{G} \rightarrow 2^{\bm{\mathcal{H}}}\times 2^{\mathcal{G}} : (\x, u) \mapsto \left(\bigtimes_{i \in \mathcal{I}} \partial \iota_{C_i}(x_i) \right) \times 
    \partial \iota_{D}^*(u). 
  \end{align}
  Then, we have the following. 
  \begin{enumerate}
    \renewcommand{\labelenumi}{\textnormal{(\roman{enumi})}}
    \item 
    $\zer(\A+\B):=\{\bxi \in \bm{\mathcal{H}}\times\mathcal{G} \mid \bm{0}_{\bm{\mathcal{H}}\times\mathcal{G}} \in \A(\bxi) + \B(\bxi)\} \bigl(=\Fix(\T_{\mathrm{FB}})\bigr) = \Fix(\TFBF) = \Fix(\T_{\FBF}^{\alpha}) \neq \varnothing$ and $\bm{\mathcal{V}}$ can be expressed as
    \begin{equation}
      \label{eq:V_fixed}
      \bm{\mathcal{V}} = \Q_{\bm{\mathcal{H}}}\bigl( \Fix(\T_{\FBF}^{\alpha}) \bigr) \Bigl( = \Q_{\bm{\mathcal{H}}}\bigl( \Fix(\T_{\mathrm{FB}}) \bigr) = \Q_{\bm{\mathcal{H}}}\bigl( \Fix(\TFBF) \bigr) \Bigr)
    \end{equation}
    with a canonical projection $\Q_{\bm{\mathcal{H}}}:\bm{\mathcal{H}} \times \mathcal{G} \rightarrow \bm{\mathcal{H}}:(\x, u)\mapsto\x$~onto~$\bm{\mathcal{H}}$. 
    \item 
    $\T_{\FBF}^{\alpha}$ is continuous and strongly attracting quasi-nonexpansive. 
    Moreover, $\T_{\FBF}^{\alpha}$ is quasi-shrinking on any bounded closed   convex set $\bm{\mathcal{C}}\subset\bm{\mathcal{H}}\times\mathcal{G}$ satisfying $\Fix(\T_{\FBF}^{\alpha})\cap\bm{\mathcal{C}} \neq \varnothing$. 
  \end{enumerate}
\end{proposition}
\textit{Proof.}
See Appendix B. 
\qed

\begin{remark}[Motivation for focusing on $\Talpha$]
  \begin{enumerate}
    \renewcommand{\labelenumi}{\textnormal{(\roman{enumi})}}
    \item (Wider applicability of $\T_{\mathrm{FBF}}$ than $\T_{\mathrm{FB}}$ for v-GNEP). 
    Fixed point theory has been utilized to design algorithms for approximating iteratively a point in the solution set $\bm{\mathcal{V}}$ of v-GNEP, see, e.g., \cite{briceno-ariasMonotoneOperatorMethods2013,Franci2020FBFforGeneralizedNash}. 
    These approaches are using the expression of $\bm{\mathcal{V}}$ in terms of the fixed point set of quasi-nonexpansive operator $\T_{\mathrm{FB}}$ or $\TFBF$ as in \eqref{eq:V_fixed}. 
    In particular, if $\A$ is $\eta (> 0)$-cocoercive\footnote{\label{fot:cocoercive}
      $\A$ is $\eta$-cocoercive if 
      \begin{equation}
      (\exists \eta > 0, \forall \bxi \in \bm{\mathcal{H}}\times\mathcal{G}, \forall \bzeta \in \bm{\mathcal{H}}\times\mathcal{G}) \ \product[\bm{\mathcal{H}}\times\mathcal{G}]{\A(\bxi) - \A(\bzeta)}{\bxi - \bzeta} \geq \eta \|\A(\bxi) - \A(\bzeta)\|^2_{\bm{\mathcal{H}}\times\mathcal{G}}.
      \end{equation}
      The above inequality ensures that $\A$ is monotone and Lipschitzian over $\bm{\mathcal{H}}\times\mathcal{G}$. 
      }, 
    then the nonexpansiveness of $\T_{\mathrm{FB}}$ in \eqref{eq:FBF} with $\gamma \in (0, 2\eta)$ is automatically guaranteed \cite[Prop. 26.1(iv)(d)]{BC2017}, and we can apply, e.g., \cite{GROETSCH1972369}, to $\T_{\mathrm{FB}}$ for approximation of a point in $\bm{\mathcal{V}}$. 
    However, the cocoercivity of $\A$ is not ensured by Assumption~\ref{asm:G_monotone_Lipschitzian}\ref{item:asm_G}\footnote{
      Assumption~\ref{asm:G_monotone_Lipschitzian} ensures the Lipschitz continuity and monotonicity of $\A$ (see Lemma~\ref{lemma:A_B_monotone}), but does not imply the cocoercivity of $\A$. 
    }, and the nonexpansiveness of $\T_{\mathrm{FB}}$ is not guaranteed automatically by Assumption~\ref{asm:G_monotone_Lipschitzian}\ref{item:asm_G}. 
    In contrast, under Assumption~\ref{asm:G_monotone_Lipschitzian}\ref{item:asm_G}, the quasi-nonexpansiveness and continuity of $\TFBF$ in \eqref{eq:FBF} are guaranteed, and we can apply, e.g., \cite[Thm. 6]{DOSTON1970}, to $\T_{\mathrm{FBF}}$ for approximation of  a point in $\bm{\mathcal{V}}$. 
    \item To be discussed in Section~\ref{subsec:algorithm}, the strongly attracting quasi-nonexpansiveness and quasi-shrinkingness\footnote{
      While the quasi-shrinkingness of $\TFBF$ can be shown by a similar discussion found in \cite[Lemma 3]{e.benenatiOptimalSelectionTracking2023}, the strongly attracting quasi-nonexpansiveness of $\TFBF$ is not ensured. 
    } of $\alpha$-averaged operator $\Talpha$ (see Proposition \ref{prop:FBF_fix}(ii)) are key ingredients in this paper for solving the upper-level v-GNEP. 
    The strongly attracting quasi-nonexpansiveness of $\Talpha$ is helpful to show the quasi-shrinkingness of $P_{\overline{B}(0;r)}\circ \T^{\alpha}_{\FBF}$, which will be used in the proposed algorithm (see Theorem 3.6). 
  \end{enumerate}
\end{remark}

Upper-level v-GNEP in \eqref{eq:upper_V} is designed as the variational inequality $\mathrm{VI}(\bm{\mathcal{V}}, \uG)$. 
Note that, to $\mathrm{VI}(\bm{\mathcal{V}}, \uG)$, we can not directly apply standard algorithms (see e.g., \cite[Exm. 26.26, Exm. 26.27]{BC2017}) using the metric projection onto $\bm{\mathcal{V}}$ because the projection onto $\bm{\mathcal{V}}$ is not available in general as a computable operator. 
To design an algorithm without requiring the metric projection onto $\bm{\mathcal{V}}$, we use a translation of $\mathrm{VI}(\bm{\mathcal{V}}, \uG)$ into a variational inequality over the fixed point set of $\T_{\FBF}^{\alpha}$ in \eqref{eq:Talpha}. 
\begin{lemma}
  \label{lem:upper_V_transformation}
  For Problem \ref{prob:H_vGNEP} under Assumption \ref{asm:G_monotone_Lipschitzian}, $\uV = \Q_{\bm{\mathcal{H}}}(\bm{\Omega})$~holds~with 
  \begin{equation}
    \label{eq:upper_VI_over_fixed_point}
    \begin{aligned}
      \bm{\Omega}\coloneqq \Bigl\{ \bxi \in \Fix(\T_{\FBF}^{\alpha}) \Bigm\vert \bigl(\forall \bzeta \in \Fix(\T_{\FBF}^{\alpha})\bigr) \ \bigl\langle \widetilde{\uG}(\bxi), \bzeta - \bxi \bigr\rangle_{\bm{\mathcal{H}}\times\mathcal{G}} \geq 0  \Bigr\} \subset \bm{\mathcal{H}}\times\mathcal{G},  
    \end{aligned}
  \end{equation}
  where (i) $\Q_{\bm{\mathcal{H}}}:\bm{\mathcal{H}}\times\mathcal{G} \rightarrow \bm{\mathcal{H}}$ is a canonical projection in Proposition \ref{prop:FBF_fix}(ii), 
  \begin{equation}
    \label{eq:extend_uG}
    (ii) \ \widetilde{\bm{\mathfrak{G}}}:\bm{\mathcal{H}}\times\mathcal{G}\rightarrow\bm{\mathcal{H}}\times\mathcal{G}:\bm{\xi}:=(\x, u) \mapsto \bigl(\bm{\mathfrak{G}}(\x), \bm{0}_{\mathcal{G}}\bigr), 
  \end{equation}
  and (iii) $\Talpha$ is defined by \eqref{eq:Talpha}. 
\end{lemma}
\textit{Proof.} 
From the definition of $\uV$ in \eqref{eq:upper_V}, 
we have
\begin{align}
  \uV &= \bigl\{\x \in \bm{\mathcal{V}} \bigm\vert
  \product[\bm{\mathcal{H}}]{\uG(\x)}{\y -\x} \geq 0 \ (\forall \y \in \bm{\mathcal{V}} \overset{\eqref{eq:V_fixed}}{=} \Q_{\bm{\mathcal{H}}}(\Fix(\T_{\FBF}^{\alpha})))\bigr\} \\
  &= \Bigl\{\x \in \Q_{\bm{\mathcal{H}}}\bigl(\Fix(\T_{\FBF}^{\alpha})\bigr) \bigm\vert (\exists u \in \mathcal{G}) \ 
  (\x, u) \in \Fix(\T_{\FBF}^{\alpha}) \ \textrm{satisfying} \\ 
  &\qquad \quad \bigl(\forall (\y, v) \in \Fix(\T_{\FBF}^{\alpha})\bigr) \ \bigl\langle \underbrace{(\uG(\x), \bm{0}_{\mathcal{G}})}_{=\widetilde{\uG}(\x, u)}, (\y -\x, v - u) \bigr\rangle_{\bm{\mathcal{H}}\times\mathcal{G}} \geq 0 \Bigr\} \\
  &= \Q_{\bm{\mathcal{H}}} \Bigl( \Bigl\{ \bxi \in \Fix(\T_{\FBF}^{\alpha}) \Bigm\vert \bigl(\forall \bzeta \in \Fix(\T_{\FBF}^{\alpha})\bigr) \ \bigl\langle \widetilde{\uG}(\bxi), \bzeta - \bxi \bigr\rangle_{\bm{\mathcal{H}}\times\mathcal{G}} \geq 0  \Bigr\} \Bigr) \\
  &= \Q_{\bm{\mathcal{H}}}(\bm{\Omega}). 
\end{align}
\qed

\subsection{Proposed Algorithm for Solving Upper-level v-GNEP}
\label{subsec:algorithm}
By Proposition \ref{prop:FBF_fix}(ii), $\T_{\FBF}^{\alpha}$ in \eqref{eq:Talpha} enjoys the quasi-shrinking condition. 
This situation and Fact \ref{fact:HSDM} encourage us to use the hybrid steepest descent method for finding a point in $\bm{\Omega}$ of \eqref{eq:upper_VI_over_fixed_point} with $\widetilde{\uG}$ in \eqref{eq:extend_uG} by generating a sequence $(\bxi_n)_{n\in\N}$ with any initial point $\bm{\xi}_0 \in \bm{\mathcal{H}}\times\mathcal{G}$:
\begin{equation}
  \label{eq:HSDM}
  (n\in\N) \ \bm{\xi}_{n+1} = \T_{\FBF}^{\alpha}(\bm{\xi}_n) - \lambda_{n+1} \widetilde{\uG}(\T_{\FBF}^{\alpha}(\bm{\xi}_n))
\end{equation}
with stepsize $(\lambda_n)_{n\in\N} \subset [0, \infty)$ satisfying\footnote{For example, $\lambda_n = 1/n$ satisfies \eqref{eq:H1_H2}.}
\begin{align}
  \label{eq:H1_H2}
  \mathrm{(H1)} \ \lim_{n\rightarrow \infty} \lambda_n = 0 \ \mathrm{and} \ \mathrm{(H2)} \ \sum_{n\in\N} \lambda_n = \infty. 
\end{align}
By using Lemma \ref{lem:upper_V_transformation}, a convergence property of the \emph{prototype algorithm} \eqref{eq:HSDM} is presented in Theorem \ref{thm:HSDM_convergence} below. 
\begin{algorithm}[t]
  \caption{Hybrid steepest descent method for Upper-level v-GNEP \eqref{eq:upper_V}}
  \label{alg:HSDM}
  \begin{algorithmic}[1]
    \State $\textbf{Input}: \gamma \in \left(0, 1/(\kappa_{\G} + \|\Li\|_{\mathrm{op}}) \right), \alpha \in (0, 1), (\lambda_n)_{n \in \N} \subset [0, \infty) \ \textrm{satisfying} \ \textrm{\eqref{eq:H1_H2}}, (\x_0, u_0) = ((x_{i, 0})_{i \in \mathcal{I}}, u_{0}) \in \bm{\mathcal{H}}\times\mathcal{G} \ \textrm{(For every} \ n \in \N, x_{i, n} \textrm{denotes the} \ i\textrm{-th component of} \ \x_n), \ r > 0 \ \textrm{satisfying} \ \Fix(\Talpha) \subset \clBall \ \textrm{(if you use the proposed algorithm \eqref{eq:replace_HSDM})}.$
    \State $\textbf{Set}: \Pi_i:\bm{\mathcal{H}}\rightarrow\mathcal{H}_i:\x\mapsto x_i. $
    \State $\textbf{for } n = 1, 2, \dots \textbf{ do}$
      \State $\quad \textrm{Forward-backward} \ \textrm{step}:$
      \State $\qquad (\forall i \in \mathcal{I}) \ y_{i, n} \gets P_{C_i} \bigl[x_{i, n} - \gamma \bigl(\nabla_i \f_i(\x_n) + \Pi_i (\Li^*u_n) \bigr) \bigr]$
      \State $\qquad w_n \gets u_n - \gamma P_{D}\bigl[(1/\gamma)u_{n} + \Li \x_n \bigr]$
      \State $\quad \textrm{Forward} \ \textrm{step}:$
      \State $\qquad (\forall i \in \mathcal{I}) \ \tilde{y}_{i, n} \gets y_{i, n} - \gamma \bigl[ \bigl(\nabla_i \f_i(\y_n) + \Pi_i (\Li^* w_n)\bigr) - \bigl(\nabla_i \f_i(\x_n) + \Pi_i (\Li^* u_n)\bigr) \bigr]$
      \State $\qquad \tilde{w}_{n} \gets w_{n} + \gamma (\Li \y_n - \Li \x_n)$
      \State $\quad \alpha\textrm{-averaging} \ \textrm{(and projection)} \ \textrm{step}:$
      \State $\quad \ \ (\x_{n+1/2}, u_{n+1}) \gets 
      \begin{cases}
        (1-\alpha)(\x_n, u_n) + \alpha(\tilde{\y}_n, \tilde{w}_n)  \ (\textrm{Prototype algorthm \eqref{eq:HSDM}}) \\
        P_{\clBall}\bigl((1-\alpha)(\x_n, u_n) + \alpha(\tilde{\y}_n, \tilde{w}_n)\bigr) \ (\textrm{Proposed algorthm \eqref{eq:replace_HSDM}})
      \end{cases}$
      \State $\quad \textrm{Steepest descent step}:$
      \State $\qquad (\forall i \in \mathcal{I}) \ x_{i, n+1} \gets x_{i, n+1/2} - \lambda_{n+1} \nabla_i \uf{i}(\x_{n+1/2})$
    \State $\textbf{end for}$
  \end{algorithmic}
\end{algorithm}

\begin{theorem}[Convergence of prototype algorithm]
  \label{thm:HSDM_convergence}
  Under Assumption \ref{asm:G_monotone_Lipschitzian}, define a nonempty closed convex set $\uV$ as in \eqref{eq:upper_V} and $\Talpha$ as in \eqref{eq:Talpha}. 
  Then, for any initial point $\bm{\xi}_0 \in \bm{\mathcal{H}}\times\mathcal{G}$, the sequence $(\bm{\xi}_n)_{n \in \N} = (\x_n, u_n)_{n\in\N}$ generated by \eqref{eq:HSDM} enjoys $\lim_{n \rightarrow \infty} d\bigl(\x_n, \uV\bigr) = 0$ provided that:
  \begin{equation}
    \label{eq:bound_cond_proto}
    \Fix(\T_{\FBF}^{\alpha})\subset\bm{\mathcal{H}}\times\mathcal{G} \ \footnote{By Proposition \ref{prop:FBF_fix}(i), $\Fix(\Talpha) \neq \varnothing$.} \ and \ (\bm{\xi}_n)_{n\in\N}\subset\bm{\mathcal{H}}\times\mathcal{G} \ are \ bounded. \tag{\ddag}
  \end{equation}
\end{theorem}
\textit{Proof.} See Appendix D. 
\qed
Unfortunately, at this stage, the condition \eqref{eq:bound_cond_proto} is not guaranteed automatically. 
In Theorem \ref{thm:modified_HSDM} below, we modify the prototype algorithm \eqref{eq:HSDM} to \eqref{eq:replace_HSDM} and show its convergence under the boundedness of $\Fix(\Talpha)$, without assuming the boundedness of the sequence generated by \eqref{eq:replace_HSDM}. 
Fortunately, the boundedness of $\Fix(\Talpha)$ is guaranteed under a mild condition (see Proposition \ref{prop:bound_cond} and Remark \ref{rem:bound_cond}). 
\begin{theorem}[Proposed algorithm: Modification of prototype algorithm and its convergence]
  \label{thm:modified_HSDM}
  Under Assumption \ref{asm:G_monotone_Lipschitzian}, suppose that there exists $r>0$ satisfying $\Fix(\Talpha) \subset \clBall\coloneqq\{\bxi\in\bm{\mathcal{H}}\times\mathcal{G} \vert \|\bxi\| \leq r\}$ for $\Talpha$ in \eqref{eq:Talpha}. 
  For any $\bm{\xi}_0 \in \bm{\mathcal{H}}\times\mathcal{G}$, generate a sequence $(\bm{\xi}_n)_{n \in \N} = (\x_n, u_n)_{n\in\N}$ by 
  \begin{equation}
    \label{eq:replace_HSDM}
    \bxi_{n+1} = \bigl( P_{\clBall} \circ \Talpha \bigr)(\bxi_n) - \lambda_{n+1} \widetilde{\uG}\Bigl(\bigl( P_{\clBall} \circ \Talpha \bigr)(\bxi_n)\Bigr),
  \end{equation}
  where 
  $\widetilde{\uG}$ is defined by \eqref{eq:extend_uG} and $(\lambda_n)_{n\in\N}$ satisfies \eqref{eq:H1_H2}. 
  Then, we have:
  \begin{enumerate}
    \renewcommand{\labelenumi}{\textnormal{(\roman{enumi})}}
    \item 
    $P_{\clBall} \circ \Talpha$ is a continuous and strongly attracting quasi-nonexpansive operator with bounded $\Fix\bigl(P_{\clBall} \circ \Talpha \bigr)=\Fix(\Talpha)\cap\clBall = \Fix(\Talpha)$. 
    Moreover, $P_{\clBall} \circ \Talpha$ is quasi-shrinking on any bounded, closed, and convex set $\bm{\mathcal{C}}\subset\bm{\mathcal{H}}\times\mathcal{G}$ satisfying $\Fix\bigl( P_{\clBall} \circ \Talpha \bigr)\cap\bm{\mathcal{C}}\neq \varnothing$. 
    \item $(\bxi_n)_{n\in\N}$ is bounded.
    \item 
    $\lim_{n \rightarrow \infty} d\bigl(\x_n, \uV \bigr) = 0$ holds true, where $\uV$ is defined in \eqref{eq:upper_V}. 
  \end{enumerate}
\end{theorem}
\textit{Proof.}
See Appendix E.
\qed

From Theorem \ref{thm:modified_HSDM}(ii), there exists a cluster point of $(\bxi_n)_{n\in\N}$ generated by \eqref{eq:replace_HSDM}, and any cluster point $\bar{\bxi} = (\bar{\x}, \bar{u})$ of $(\bxi_n)_{n\in\N} = (\x_n, u_n)_{n\in\N}$ belongs to $\bm{\Omega}$ in \eqref{eq:upper_VI_over_fixed_point} (see Remark \ref{rem:HSDM}(i)). 
In addition, we have $\bar{\x} \in \uV$ by Lemma \ref{lem:upper_V_transformation}. 
Algorithm \ref{alg:HSDM} illustrates concrete expressions of the prototype algorithm \eqref{eq:HSDM} and of the proposed algorithm \eqref{eq:replace_HSDM} (the two algorithms differ only in Line 11). 
We now turn to the boundedness of $\Fix(\Talpha)$ below. 
\begin{proposition}[A sufficient condition to ensure the boundedness of $\Fix(\Talpha)$]
  \label{prop:bound_cond}
  Under Assumption \ref{asm:G_monotone_Lipschitzian}, suppose that:
  \begin{enumerate}[label=\textnormal{(\alph*)}]
    \item \label{item:bound_cond2} $\bm{\mathcal{V}}$ in \eqref{eq:VE} is bounded (this is guaranteed automatically when $\bm{\mathfrak{C}}$ in \eqref{eq:fracture_C} is bounded). 
    \item \label{item:bound_cond3} Linear operator $\Li:\bm{\mathcal{H}}\to\mathcal{G}$ has full column rank, i.e., $\operatorname{ker}(\Li^*) = \{0_{\mathcal{G}}\}$. 
    \item \label{item:bound_cond4} For every $\x \in \bm{C}\cap\Li^{-1}(D) = \bm{\mathfrak{C}}$, $\partial\iota_{\C}(\x) \cap - \partial\iota_{\Li^{-1}(D)}(\x) = \{\bm{0}_{\bm{\mathcal{H}}}\}$, i.e., 
    \begin{equation}
      \label{eq:cond_sumrule_R}
      \begin{aligned}
        &\bigl(\forall \x \in \bm{C}\cap\Li^{-1}(D), \forall \bm{v}_1 \in \partial\iota_{\C}(\x), \forall \bm{v}_2 \in \partial\iota_{\Li^{-1}(D)}(\x)\bigr) \\ 
        &\qquad\qquad\qquad\qquad\qquad\qquad
        \vv_1 + \vv_2 = \bm{0}_{\bm{\mathcal{H}}} 
        \Rightarrow 
        \vv_1 = \vv_2 = \bm{0}_{\bm{\mathcal{H}}}.
      \end{aligned}
  \end{equation}
  \end{enumerate}
  Then, for $\Talpha$ in \eqref{eq:Talpha}, $\Fix(\Talpha)$ is bounded. 
\end{proposition}
\textit{Proof.}
  See Appendix F. 
\qed
\begin{remark}
  \label{rem:bound_cond} \ 
  \begin{enumerate}[label=(\roman*)]
    \item (Reasonability of Proposition \ref{prop:bound_cond}\ref{item:bound_cond2}). The condition in Proposition \ref{prop:bound_cond}\ref{item:bound_cond2} is one of sufficient conditions to ensure $\uV \neq \varnothing$ (see Fact \ref{fact:vGNE_GNEP}(i)). 
    \item (Relation~between~Proposition~\ref{prop:bound_cond}(b)-(c)~and~qualification~condition). 
    \begin{enumerate}[label=\textnormal{(\roman{enumi}-\arabic*)}]
      \item \label{item:rem_sumrule} Condition \eqref{eq:cond_sumrule_R} is well known as a qualification condition (see, e.g., \cite[Cor. 10.9]{RW1998}) for $\partial (\iota_{\C} + \iota_{\Li^{-1}(D)}) = \partial \iota_{\C} + \partial \iota_{\Li^{-1}(D)}$. 
    \item \label{item:rem_chainrule} Under the condition in Proposition \ref{prop:bound_cond}\ref{item:bound_cond3}, we have \linebreak$\mathrm{ri}(\dom(\iota_D) - \mathrm{ran}(\Li)) = \mathrm{ri}(D - \mathcal{G}) = \operatorname{ri}(\mathcal{G}) = \mathcal{G}\ni 0$.
    Then, by Fact \ref{fact:subdifferential}\ref{fact:chain_rule}, we have $\partial (\iota_D \circ \Li) = \Li^* \circ \partial \iota_D \circ \Li$. 
    \item By combining 
    $\iota_{\Li^{-1}(D)} = \iota_D\circ\Li$
    with \ref{item:rem_sumrule} and \ref{item:rem_chainrule}, Assumption \ref{asm:G_monotone_Lipschitzian}\ref{item:asm_sum_chain} is guaranteed automatically under the conditions in Proposition \ref{prop:bound_cond}\ref{item:bound_cond3} and \ref{item:bound_cond4}. 
    \end{enumerate}
  \end{enumerate}
\end{remark}
The next result with $(S_1, S_2) \coloneqq \bigl(\C, \Li^{-1}(D)\bigr)$ or $(S_1, S_2) \coloneqq \bigl(\Li^{-1}(D), \C\bigr)$ gives sufficient conditions to ensure the condition \eqref{eq:cond_sumrule_R}. 
\begin{proposition}
  \label{prop:cond_sumrule_R-BC}
  Let $\mathcal{H}$ be a finite-dimensional real Hilbert space, and $S_1 \subset \mathcal{H}$ and $S_2 \subset \mathcal{H}$ be nonempty closed convex sets satisfying $S_1 \cap S_2 \neq \varnothing$. 
  Suppose that \ref{item:QC_BC_1} or \ref{item:QC_BC_2} below holds. 
  \begin{enumerate}[label=\textnormal{(\alph*)}]
    \item \label{item:QC_BC_1} $\operatorname{int}(S_1) \cap S_2 \neq \varnothing$.
    \item \label{item:QC_BC_2} All of the following hold.
    \begin{enumerate}[label=\textnormal{(\alph{enumi}-\arabic*)}]
      \item \label{subitem:QC_ri_C1} $\operatorname{ri}(S_1) \cap S_2 \neq \varnothing$.
      \item \label{subitem:QC_ri_C2} $S_1 \cap \operatorname{ri}(S_2) \neq \varnothing$.
      \item \label{subitem:QC_ri_aff} $\operatorname{aff}(S_1) + \operatorname{aff}(S_2) = \mathcal{H}$. 
    \end{enumerate}
  \end{enumerate}
  Then we have 
  \begin{equation}
    \label{eq:QC_R}
    (\forall x \in S_1 \cap S_2) \ \partial \iota_{S_1}(x) \cap - \partial \iota_{S_2}(x) = \{0_{\mathcal{H}}\}.
  \end{equation}
\end{proposition}
\textit{Proof.}
See Appendix G.
\qed

\begin{remark}
  \label{rem:qualification_cond} \ 
  \begin{enumerate}[label=\textnormal{(\roman*)}]
    \item (Interpretation of Proposition \ref{prop:cond_sumrule_R-BC}). The condition in Proposition \ref{prop:cond_sumrule_R-BC}\ref{item:QC_BC_1} is well known as a qualification condition for $\partial (\iota_{S_1} + \iota_{S_2}) = \partial \iota_{S_1} + \partial \iota_{S_2}$ \cite[Cor. 16.48(ii)]{BC2017}. 
    In addition, the conditions in Proposition \ref{prop:cond_sumrule_R-BC}\ref{subitem:QC_ri_C1} and \ref{subitem:QC_ri_C2} are satisfied if another qualification condition $\operatorname{ri}(S_1) \cap \operatorname{ri}(S_2) \neq \varnothing$ \cite[Cor. 16.48(iv)]{BC2017} holds true. 
    \item (Note on sufficient conditions for \eqref{eq:QC_R}). The condition $0_{\mathcal{H}} \in \operatorname{ri}(S_1 - S_2)$ has been reported (e.g., \cite[Remark 1]{huAnalysisSemismoothNewtonType2025}) as a sufficient condition to ensure \eqref{eq:QC_R} without any proof; however, this is not always valid. 
    To examine this, for $\mathcal{H}\coloneqq \R^2$, let $S_1 = S_2 = \{(s, 0) \in \R^2 \mid s \in \R \}$. 
    Then, we have $\operatorname{ri}(S_1 - S_2) = S_1 \ni 0_{\R^2}$. 
    Nevertheless, we have $\partial \iota_{S_1}(x) = \partial \iota_{S_2}(x) = \{(0, t) \in \R^2 \mid t \in \R\}$ for every $x \in S_1 = S_2$, which implies that $\partial \iota_{S_1}(x) \cap - \partial \iota_{S_2}(x) = \{(0, t) \in \R^2 \mid t \in \R\} \neq \{0_{\R^2}\}$.
  \end{enumerate}
\end{remark}

\subsection{Equilibrium Selection via Cycle}
\label{subsec:cycle}
Consider the lower-level non-cooperative game formulated in the form of a special case of Problem \ref{prob:GNEP} where $m\geq 2$ and $\mathcal{H}_i = \mathcal{H} \ (i \in \mathcal{I})$ with a common finite-dimensional real Hilbert space $\mathcal{H}$. 
Assume that $\bm{\mathcal{V}}$ in \eqref{eq:VE} is a nonempty closed convex set. 
As an equilibrium selection problem motivated by cycles (see Section \ref{subsec:preriminary_cycle}), we consider finding a \underline{cycle over $\bm{\mathcal{V}}$}, in the following sense (see Definition \ref{def:cycle} and Lemma \ref{lem:cycle_V_cycle}): 
\begin{equation}
  \label{eq:selection_cycle}
  \begin{aligned}
    &\find \bar{\x} \in \mathrm{cyc}^{\langle \mathrm{u} \rangle}(\bm{\mathcal{V}}) \coloneqq \bigl\{ (\bar{x}_1, \dots, \bar{x}_m) \in \bm{\mathcal{V}} \bigm\vert  \\ 
    &\bar{x}_2 = P_{\mathcal{V}_2(\bar{\x}_{\smallsetminus 2})}(\bar{x}_1), \dots, \bar{x}_{m} = P_{\mathcal{V}_{m}(\bar{\x}_{\smallsetminus m})}(\bar{x}_{m-1}), \bar{x}_1 = P_{\mathcal{V}_1(\bar{\x}_{\smallsetminus 1})}(\bar{x}_{m}) \bigr\}, 
  \end{aligned}
\end{equation}
where $\mathcal{V}_i(\x_{\smallsetminus i}) \coloneqq \{x_i\in\mathcal{H} \vert (x_i; \x_{\smallsetminus i}) \in \bm{\mathcal{V}}\} \ (i\in\mathcal{I})$ are nonempty closed convex sets for every $\x\in\bm{\mathcal{V}}$. 
Indeed, cycles over $\bm{\mathcal{V}}$ formulated as solutions of Problem \eqref{eq:selection_cycle} are an extension of cycles in Definition \ref{def:cycle}. 

\begin{lemma}\label{lem:cycle_V_cycle}
  Consider $\mathrm{cyc}^{\langle \mathrm{u} \rangle}(\bm{\mathcal{V}})$ in \eqref{eq:selection_cycle} 
  for a special case where $\bm{\mathcal{V}}$ in \eqref{eq:VE} can be expressed as $\bm{\mathcal{V}}=\bigtimes_{i\in\mathcal{I}} K_i$ in terms of nonempty closed convex sets $K_i \subset \mathcal{H} \ (i\in\mathcal{I})$. 
  Then $\mathrm{cyc}^{\langle \mathrm{u} \rangle}(\bm{\mathcal{V}}) = \mathrm{cyc}(K_1, \dots, K_m)$ (see also Definition \ref{def:cycle}). 
\end{lemma}
\textit{Proof.} 
Clear from $\mathcal{V}_i(\x_{\smallsetminus i}) = K_i \ (i\in\mathcal{I})$. 
\qed

\begin{remark}[Interpretation of Problem \eqref{eq:selection_cycle}]
  By analogy with the argument in footnote \footref{fot:cycle_game_eq} in the end of Section \ref{subsec:preriminary_cycle}, $\mathrm{cyc}^{\langle \mathrm{u} \rangle}(\bm{\mathcal{V}})$ in \eqref{eq:selection_cycle} can be expressed as the following GNEP:
    \begin{equation}
      \label{eq:HGNEP_cycle}
      \begin{aligned}
        &\mathrm{cyc}^{\langle \mathrm{u} \rangle}(\bm{\mathcal{V}}) = \SGNE{\mathcal{H}^m}{\mathcal{G}}{\uf{i}}{\mathcal{H}}{\Id}{\bm{\mathcal{V}}} \\
        &= \Bigl\{ (\bar{x}_1, \dots, \bar{x}_m) \in \mathcal{H}^m \Bigm\vert \bigl(\forall i \in \mathcal{I}\bigr) \ \bar{x}_i \in \argmin[x_i\in\mathcal{H}] \ \uf{i}(x_i;\bar{\x}_{\smallsetminus i}) \ \st (x_i; \bar{\x}_{\smallsetminus i}) \in \bm{\mathcal{V}} \Bigr\}
      \end{aligned}
    \end{equation}
    with the convention $0 = m$, where 
    \begin{equation}
      \label{eq:uf_cycle}
      (i\in\mathcal{I}) \ \uf{i}: \mathcal{H}^m \to \R: \x = (x_1, \dots, x_m) \mapsto \frac{1}{2}\|x_i - x_{i-1}\|^2. 
    \end{equation}
    In other words, Problem \eqref{eq:selection_cycle} requires to find an $\bar{\x}=(\bar{x}_1, \dots, \bar{x}_m)$ such that player $i$'s strategy $\bar{x}_i$ is located in $\mathcal{V}(\bar{\x}_{\smallsetminus i})$ closest to player $(i-1)$'s strategy $\bar{x}_{i-1}$. 
    We believe that such a selection by Problem \eqref{eq:selection_cycle} has potential being used in distributed control of some multiagent systems where closeness among players' strategies is desirable (see, e.g., \cite[Sec. VI.A]{yeDistributedNashEquilibrium2023}). 
\end{remark}

Lemma \ref{lem:cycle_V_cycle} implies that the POCS algorithm (see Remark \ref{rem:POCS}) is applicable directly to Problem \eqref{eq:selection_cycle} if (a) $\bm{\mathcal{V}}$ in \eqref{eq:VE} can be expressed as $\bm{\mathcal{V}}=\bigtimes_{i\in\mathcal{I}} K_i$ in terms of nonempty closed convex sets $K_i \subset \mathcal{H} \ (i\in\mathcal{I})$, and (b) all $P_{K_i} \ (i\in\mathcal{I})$ are available.
However, Problem \eqref{eq:selection_cycle} in general cases is challenging because $P_{\mathcal{V}_i(\x_{\smallsetminus i})} \ (i\in \mathcal{I})$ are not available (see, e.g., Example \ref{ex:cycle_implicit}) as computable operators. 
We start with an expression of Problem \eqref{eq:selection_cycle} in terms of the GNEP \eqref{eq:HGNEP_cycle} and consider an application of the proposed algorithm \eqref{eq:replace_HSDM} to the \underline{v-GNEP for $\sSGNE{\mathcal{H}^m}{\mathcal{G}}$ in \eqref{eq:HGNEP_cycle}} as an upper-level v-GNEP over $\bm{\mathcal{V}}$. 
\begin{problem}[Upper-level v-GNEP for finding a cycle over $\bm{\mathcal{V}}$]
  \label{prob:cycle_selection}
  Under the setting of the lower-level non-cooperative game formulated in the form of a special case of Problem \ref{prob:GNEP}, where $m\geq 2$ and $\mathcal{H}_i = \mathcal{H} \ (i \in \mathcal{I})$ with a common finite-dimensional real Hilbert space $\mathcal{H}$, assume that $\bm{\mathcal{V}}$ in \eqref{eq:VE} is a nonempty closed convex set\footnote{
    See Fact \ref{fact:vGNE_GNEP}(i) for its sufficient condition. 
  }. 
  The upper-level v-GNEP for finding a cycle over $\bm{\mathcal{V}}$ is given as v-GNEP for $\sSGNE{\mathcal{H}^m}{\mathcal{G}}$ in \eqref{eq:HGNEP_cycle} with $\uf{i} \ (i\in\mathcal{I})$ in \eqref{eq:uf_cycle}:
  \begin{equation}
    \label{eq:HGNEP_cycle_V}
    \find \x^{\star} \in \uV_{\mathrm{cyc}} \coloneqq \bigl\{\x \in \bm{\mathcal{V}} \bigm\vert
    \product[\bm{\mathcal{H}}]{\uG_{\mathrm{cyc}}(\x)}{\y -\x} \geq 0 \ (\forall \y \in \bm{\mathcal{V}})\bigr\}, 
  \end{equation}
  where $\uV_{\mathrm{cyc}}\neq\varnothing$ is assumed\footnote{
    The boundedness of $\bm{\mathcal{V}}$ is one of sufficient conditions (see Fact \ref{fact:vGNE_GNEP}(i)).
  } and, 
  \begin{equation}
    \label{eq:G_cycle}
    \uG_{\mathrm{cyc}}:\bm{\mathcal{H}}\to\bm{\mathcal{H}}:\x  \mapsto \bigl(\nabla_1 \uf{1}(\x), \dots, \nabla_m\uf{m}(\x)\bigr).
  \end{equation}
\end{problem}

Since Problem \ref{prob:cycle_selection} is a special instance of Problem \ref{prob:H_vGNEP}, Fact \ref{fact:vGNE_GNEP}(ii) yields 
\begin{equation}
  \label{eq:cycle_V_cycle}
  \uV_{\mathrm{cyc}} \subset \SGNE{\mathcal{H}^m}{\mathcal{G}}{\uf{i}}{\mathcal{H}}{\Id}{\bm{\mathcal{V}}} \bigl(= \mathrm{cyc}^{\langle \mathrm{u} \rangle}(\bm{\mathcal{V}})\bigr). 
\end{equation}

Under Assumption \ref{asm:G_monotone_Lipschitzian}\ref{item:asm_sum_chain} and \ref{item:asm_G}, we can apply the proposed algorithm \eqref{eq:replace_HSDM} to Problem \eqref{eq:HGNEP_cycle_V} (see Corollary \ref{cor:cycle_selection} below) because Assumption \ref{asm:G_monotone_Lipschitzian}\ref{item:asm_frakG} with $\uG \coloneqq \uG_{\mathrm{cyc}}$ in \eqref{eq:G_cycle} is guaranteed automatically as follows (regarding the boundedness of $\Fix(\Talpha)$, see Proposition \ref{prop:bound_cond}). 
\begin{proposition}\label{prop:uG_cycle}
  The operator $\uG_{\mathrm{cyc}}$ defined as in \eqref{eq:G_cycle}, with $\uf{i} \ (i\in\mathcal{I})$ in \eqref{eq:uf_cycle}, is paramonotone and Lipschitzian over $\bm{\mathcal{H}}$. 
\end{proposition}
\textit{Proof.} 
By $\nabla_i\uf{i}(x_1, \dots, x_m) = x_i - x_{i-1}$, $\uG_{\mathrm{cyc}}$ can be expressed as $\uG_{\mathrm{cyc}} = \Id - \bm{\mathcal{R}}$, where the circular right-shift operator $\bm{\mathcal{R}}:\bm{\mathcal{H}}\to\bm{\mathcal{H}}:(x_1, \dots, x_m)\mapsto(x_m, x_1, \dots, x_{m-1})$ is nonexpansive. 
Therefore, $\uG_{\mathrm{cyc}}$ is Lipschitzian, and paramonotone over $\bm{\mathcal{H}}$ by \cite[Exm. 22.9]{BC2017}. 
\qed

\begin{corollary}
  \label{cor:cycle_selection}
  For Problem \ref{prob:cycle_selection} under Assumption \ref{asm:G_monotone_Lipschitzian}\ref{item:asm_sum_chain} and \ref{item:asm_G}, assume that $\Fix(\Talpha)$ is bounded\footnote{
    For its sufficient condition, see Proposition~\ref{prop:bound_cond}. 
  } for $\Talpha$ in \eqref{eq:Talpha}. 
  Generate a sequence $(\bm{\xi}_n)_{n \in \N} = (\x_n, u_n)_{n\in\N} \subset \bm{\mathcal{H}}\times\mathcal{G}$ by \eqref{eq:replace_HSDM} with $\widetilde{\uG}:\bm{\mathcal{H}}\times\mathcal{G} \to \bm{\mathcal{H}}\times\mathcal{G} : (\x, u) \mapsto (\uG_{\mathrm{cyc}}(\x), \bm{0}_{\mathcal{G}})$ and a sufficiently large $r >0 $ such that $\Fix(\Talpha) \subset \clBall$. 
  Then, for $\uV_{\mathrm{cyc}}$ in \eqref{eq:HGNEP_cycle_V}, $(\bm{\xi}_n)_{n \in \N} = (\x_n, u_n)_{n\in\N}$ enjoys $\lim_{n\to\infty} d\bigl(\x_n, \uV_{\mathrm{cyc}}\bigr) = 0$. 
\end{corollary}
\textit{Proof.}  
By Proposition \ref{prop:uG_cycle}, we can apply Theorem \ref{thm:modified_HSDM} to Problem \ref{prob:cycle_selection}. 
Then, we have $\lim_{n\to\infty} d\bigl(\x_n, \uV_{\mathrm{cyc}}\bigr) = 0$ from Theorem~\ref{thm:modified_HSDM}(iii). 
\qed

\begin{example}[Cycles associated with solution sets of smooth convex optimization problems]
  \label{ex:cycle_implicit}
  Consider the following problem:
  \begin{equation}
    \label{eq:cycle_implicit}
    \find \bar{\x} \in \mathrm{cyc}(K_1, \dots, K_m), 
  \end{equation}
  where the ordered family of nonempty closed convex sets $(K_1, \dots, K_m)$ is given implicitly, with differentiable convex functions $h_i:\mathcal{H}\to\R \ (i\in\mathcal{I})$, by
  \begin{equation}
    \label{eq:cycle_implicit_optimization_def}
    (K_1, \dots, K_m) \coloneqq (\mathrm{argmin}_{x_1\in\mathcal{H}} \ h_1(x_1),\dots, \mathrm{argmin}_{x_m\in\mathcal{H}} \ h_m(x_m)).
  \end{equation}
  For naive approaches to Problem \eqref{eq:cycle_implicit}, we remark below:
  \begin{enumerate}
    \renewcommand{\labelenumi}{\textnormal{(\roman{enumi})}}
    \item 
    POCS algorithm (see Remark \ref{rem:POCS}) can not be applied directly to Problem \eqref{eq:cycle_implicit} because computabilities of $P_{K_i} \ (i\in\mathcal{I})$ for \eqref{eq:cycle_implicit_optimization_def} are questionable. 
    \item Reducability of Problem \eqref{eq:cycle_implicit} to a hierarchical convex optimization problem \eqref{eq:upper_V_common} is questionable by Fact \ref{fact:no_variational_cycle}. 
    \item 
    $\mathrm{cyc}(K_1, \dots, K_m)$ in \eqref{eq:cycle_implicit} can be expressed\footnote{
      The solution set $\bm{\mathcal{V}}$ of $\VI(\mathcal{H}^m, \G)$ can be expressed as 
    \begin{equation}
      \label{eq:cycle_implicit_V}
      \begin{aligned}
        \bm{\mathcal{V}} 
        &= \Bigl\{ \x = (x_1, \dots, x_m) \in \mathcal{H}^m \vert (\forall \y \in \mathcal{H}^m) \ \bigl\langle (\nabla h_1(x_1), \dots, \nabla h_m(x_m)), \y - \x \bigr\rangle_{\mathcal{H}^m} \geq 0\Bigr\} \\
        &= \{ (x_1, \dots, x_m) \in \mathcal{H}^m \vert (\nabla h_1(x_1), \dots, \nabla h_m(x_m)) = \bm{0}_{\mathcal{H}^m}\} \\
        &= \{ (x_1, \dots, x_m) \in \mathcal{H}^m \vert (\forall i\in\mathcal{I}) \ \nabla h_i(x_i) = \bm{0}_{\mathcal{H}}\} \\
        &= \{ (x_1, \dots, x_m) \in \mathcal{H}^m \vert (\forall i\in\mathcal{I}) \ x_i \in \mathrm{argmin}_{z_i\in\mathcal{H}}h_i(z_i)\}
        = \textstyle\bigtimes_{i\in\mathcal{I}} K_i.
      \end{aligned}
    \end{equation}
    By $\bm{\mathcal{V}} = \bigtimes_{i\in\mathcal{I}}K_i$ and Lemma \ref{lem:cycle_V_cycle}, $\mathrm{cyc}(K_1, \dots, K_m)=\mathrm{cyc}^{\langle \mathrm{u} \rangle}(\bm{\mathcal{V}}) \overset{\eqref{eq:HGNEP_cycle}}{=} \SGNE{\mathcal{H}^m}{\mathcal{G}}{\uf{i}}{\mathcal{H}}{\Id}{\bm{\mathcal{V}}}$. 
    By Fact \ref{fact:vGNE_GNEP}(iii) with $\bm{\mathcal{V}} = \bigtimes_{i\in\mathcal{I}}K_i$, we have $\SGNE{\mathcal{H}^m}{\mathcal{G}}{\uf{i}}{\mathcal{H}}{\Id}{\bm{\mathcal{V}}} = \uV_{\mathrm{cyc}}$ in \eqref{eq:HGNEP_cycle_V}. 
    Then, we have $\mathrm{cyc}(K_1, \dots, K_m) = \uV_{\mathrm{cyc}}$.
    } as $\uV_{\mathrm{cyc}}$ in 
    Problem \ref{prob:cycle_selection} by letting $\bm{\mathcal{V}}$ be the solution set of v-GNEP $\VI(\mathcal{H}^m, \G)$ for \linebreak$\SGNE{\mathcal{H}^m}{\mathcal{H}^m}{\bm{h}_i}{\mathcal{H}}{\Id}{\mathcal{H}^m}$ with $\bm{h}_i :\mathcal{H}^m\to\R:\x \mapsto h_i(x_i) \ (i\in\mathcal{I})$ and $\G:\mathcal{H}^m\to\mathcal{H}^m:\x\mapsto(\nabla_1 \bm{h}_1(\x), \dots, \nabla_m \bm{h}_m(\x)) = (\nabla h_1(x_1), \dots, \nabla h_m(x_m))$. 
    Therefore, we can solve Problem \eqref{eq:cycle_implicit} by applying the proposed algorithm \eqref{eq:replace_HSDM} (see Corollary \ref{cor:cycle_selection}) if (a) each $\nabla h_i:\mathcal{H}\to\mathcal{H} \ (i\in\mathcal{I})$ is Lipschitzian over $\mathcal{H}$ and (b) $\Fix(\Talpha)$ is bounded, where $\Talpha$ is defined as in \eqref{eq:Talpha}. 
  \end{enumerate}
\end{example}
Example \ref{ex:cycle_implicit} tells us the remarkable expressive ability of Problem \ref{prob:H_vGNEP}. 
Moreover, Theorem \ref{thm:modified_HSDM} tells us that we can apply Algorithm \eqref{eq:replace_HSDM} to Problem \ref{prob:H_vGNEP} in a unified way. 

%% file: docs/Application.tex
\section{Numerical Experiments}
\label{sec:experiments}
To illustrate Problem \ref{prob:H_vGNEP} and the proposed algorithm \eqref{eq:replace_HSDM}, we present numerical experiments in two scenarios: (i) cycles associated with solution sets of smooth convex optimization problems (see Example \ref{ex:cycle_implicit}), and (ii) equilibrium selections from the solution set of a v-GNEP for a linearly coupled game (see, e.g., \cite{g.belgioiosoDistributedProximalPoint2019,e.benenatiOptimalSelectionGeneralized06}). 
\vspace{-1em}
\subsection{Cycles Associated with Solution Sets of Smooth Convex Optimization Problems}
To verify whether Algorithm \eqref{eq:replace_HSDM} can approximate iteratively a solution of the problem \eqref{eq:cycle_implicit} in Example \ref{ex:cycle_implicit}, we used an explicit setting\footnote{$h_i \ (i\in\mathcal{I})$ are chosen to enjoy that projections onto $K_i (= \mathrm{argmin}_{x_i\in\mathcal{H}} \ h_i(x_i))$ are available as computable operators. }: $m \coloneqq 6$, $\mathcal{H} \coloneqq \R^3$, $h_i \coloneqq \frac{1}{2} d(\cdot, K_i)^2 \ (i \in \mathcal{I})$, and $K_i \coloneqq \bigtimes_{k \in \{1, \dots, 3\}}[b_{i, k}^{\mathrm{low}}, b_{i, k}^{\mathrm{up}}] \subset [0, 100]^3 \ (i \in \mathcal{I})$ such that $(\forall i \in \mathcal{I}) \ K_i \neq \varnothing$ and $\bigcap_{i\in\mathcal{I}}K_i = \varnothing$, where $h_i \ (i\in\mathcal{I})$ are convex differentiable functions with $\nabla h_i = \Id - P_{K_i} \ (i\in\mathcal{I})$, which are Lipschitzian over $\R^3$. 
In Algorithm \eqref{eq:replace_HSDM}, we employed $(\gamma, \alpha, r) = (0.2, 0.5, 10^{15})$, $\lambda_n = 1/n$, and $\bxi_0 = \bm{0}_{\mathcal{H}^m \times \mathcal{H}^m}$. 

\begin{figure}[t]
  \centering
  \begin{minipage}[b]{0.45\columnwidth}
      \centering
      \includegraphics[clip, width=\columnwidth]{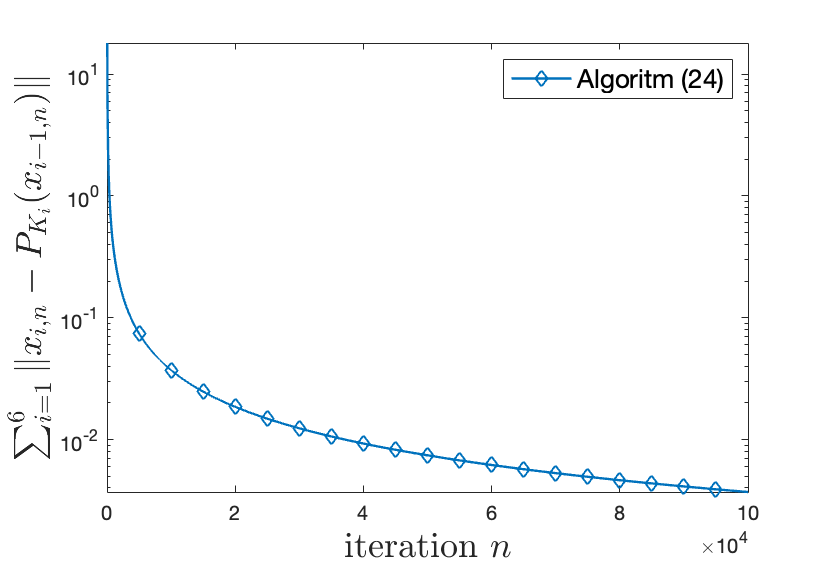}
      \caption{Cycle achievement level at each iteration.}
      \label{fig:ex1_residual}
  \end{minipage}
  \hspace{0.01\columnwidth}
  \begin{minipage}[b]{0.52\columnwidth}
      \centering
      \includegraphics[clip, width=0.85\columnwidth]{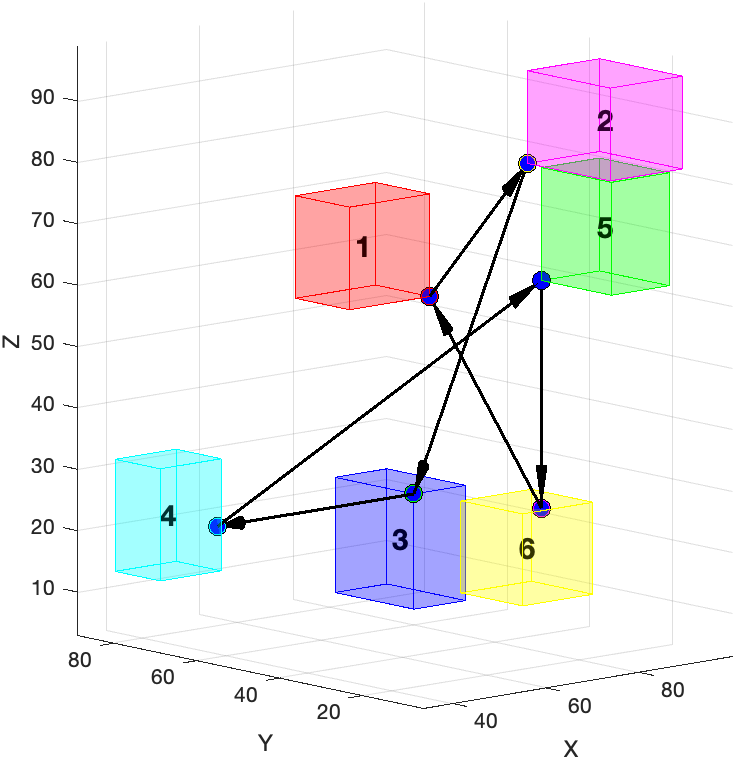}
      \caption{Visualization of $(K_1, \dots ,K_6)$ and an approximated cycle achieved by Algorithm \eqref{eq:replace_HSDM}.}
      \label{fig:ex1_3D}
  \end{minipage}
\end{figure}

We employed $\sum_{i=1}^{6}\|x_{i, n} - P_{K_i}(x_{i-1, n})\|$ to evaluate the achieving level at $\x_n= (x_{1, n}, \dots, x_{6, n})$, generated by Algorithm \eqref{eq:replace_HSDM}, toward $\mathrm{cyc}(K_1, K_2, ..., K_6)$ motivated by the fact $\sum_{i=1}^6 \|{x_{i,n}-P_{K_i}(x_{i-1,n})}\|=0 \Leftrightarrow \x_{n} \in \mathrm{cyc}(K_1, K_2, ..., K_6)$. 
Fig. \ref{fig:ex1_residual} shows the value of $\sum_{i=1}^{6}\|x_{i, n} - P_{K_i}(x_{i-1, n})\|$. 
From Fig. \ref{fig:ex1_residual}, we can see that the value of $\sum_{i=1}^{6}\|x_{i, n} - P_{K_i}(x_{i-1, n})\|$ approaches zero as $n$ increases. 
Fig. \ref{fig:ex1_3D} visualizes $K_i \ (i\in\{1, \dots, 6\})$ and $(x_1, ..., x_6)$ obtained by Algorithm \eqref{eq:replace_HSDM} with $10^5$ iterations. 
From Fig. \ref{fig:ex1_residual} and Fig. \ref{fig:ex1_3D}, we can see that Algorithm \eqref{eq:replace_HSDM} approximates iteratively a cycle associated with $(K_1, \dots, K_6)$. 

We remark that as mentioned in Example \ref{ex:cycle_implicit}(iii), Algorithm \eqref{eq:replace_HSDM} is applicable to \eqref{eq:cycle_implicit} if $\nabla h_i \ (i\in\mathcal{I})$ are available as computable operators [see, e.g., Section \ref{subsec:ex2}(Case 2) for more general case]. 

\subsection{Equilibrium Selections from Solution Set of v-GNEP for Linearly Coupled Game}
\label{subsec:ex2}

A linearly coupled game (e.g., \cite{g.belgioiosoDistributedProximalPoint2019,e.benenatiOptimalSelectionGeneralized06}) is given as an instance of Problem \ref{prob:GNEP}:
\begin{align}
  \label{eq:Cournot-Nash}
  \find \bar{\x} \in \SGNE{\R^{Mm}}{\R^M}{\f_i}{C_i}{\Li}{D},
\end{align}
where for every $i \in \mathcal{I}$,\footnotetext[12]{$y \leq c$ means that $y_j \leq c_j$ for each component $j\in\{1, \dots, M\}$.}
\begin{align}
  &\f_i(\x) := \biggl(\sum_{k\in \mathcal{I}} \bm{W}_k x_k - p \biggr)^Tx_i, \ 
  C_i := \bigtimes_{j \in \{1, \dots, M\}}[b_{i, j}^{\mathrm{low}}, b_{i, j}^{\mathrm{up}}] (\neq \varnothing) \subset \R^M \\
  &\Li:\R^{Mm}\rightarrow\R^M:\x:=(x_1, \dots, x_m) \mapsto \sum_{i\in\mathcal{I}} x_i, \ \mathrm{and}\footnotemark[12] \  D := \{y\in\R^M\mid y \leq c\}
\end{align}
\noindent with $p \in \R^M$, nonnegative diagonal matrices $\bm{W}_k\in\R^{M\times M} \ (k\in\mathcal{I})$, $b_{i, j}^{\mathrm{low}} \leq b_{i, j}^{\mathrm{up}}$, and a component-wise upper bound $c=(c_1, \dots, c_M)\in\R^M_{++}$ of $\Li(\x)\in\R^M$. 
To ensure $\bigtimes_{i\in\mathcal{I}} C_i \cap \Li^{-1}(D) \neq \varnothing$, we used $c_j > \sum_{i\in\mathcal{I}} b_{i, j}^{\mathrm{low}}$ for $j\in\{1, \dots, M\}$. 
In this case, (i) $\bm{\mathcal{V}} \neq \varnothing$, (ii) Assumption \ref{asm:G_monotone_Lipschitzian}\ref{item:asm_sum_chain} and \ref{item:asm_G}, and (iii) the boundedness of $\Fix(\Talpha)$ are satisfied. 
In numerical experiments, we used\footnote{We set the parameters in Problem \eqref{eq:Cournot-Nash} along the setting found in \cite{e.benenatiOptimalSelectionGeneralized06}, where an equilibrium selection over $\bm{\mathcal{V}}$ is considered as a convex optimization based on a hybrid steepest descent method \cite{oguraNonstrictlyConvexMinimization2003}.} $(m, M) = (6, 3)$, $c = (120, 120, 120)^T$, 
$p \in [0, 10]^3$, $b_{i, j}^{\mathrm{low}} \in [-1, 1]$, $b_{i, j}^{\mathrm{up}} = 100$ $(i \in \mathcal{I}, j \in \{1, \dots, M\})$, and $\bm{W}_k \ (k\in\mathcal{I})$ whose diagonal entries were chosen from $[0, 1]$. 

For equilibrium selection from the solution set $\bm{\mathcal{V}}$ of v-GNEP for \linebreak $\SGNE{\R^{Mm}}{\R^M}{\f_i}{C_i}{\Li}{D}$, we examined following two settings (Case 1 and Case 2) of Problem \ref{prob:H_vGNEP} (in each case, $\uV$ in \eqref{eq:upper_V} is nonempty). 
\begin{enumerate}[label=\textnormal{(\roman*)}, leftmargin=40pt]
  \renewcommand{\labelenumi}{\textnormal{(\roman{enumi})}}
  \item[(Case 1)]
  Motivated by \cite[Sec. VI.A]{yeDistributedNashEquilibrium2023}, the upper-level v-GNEP in \eqref{eq:upper_V}~is~given~with 
  \begin{equation}
    \label{eq:evaluation}
    (i\in\mathcal{I}) \ \uf{i}(\x) := \frac{1}{2} \biggl( \|x_i - t_i\|^2 + \sum_{j \in \mathcal{I} \setminus \{i\}} \|x_i - x_j\|^2 \biggr),
  \end{equation}
  where each target vector $t_i \in \R^M \ (i\in\mathcal{I})$ is chosen by player $i$. 
  Case~1 satisfies Assumption~\ref{asm:G_monotone_Lipschitzian}\ref{item:asm_frakG}. 
  In this experiment, $t_i$ is chosen from $C_i$. 
  \item[(Case 2)]
  Motivated by the discussion in Section \ref{subsec:cycle}, the upper-level v-GNEP is given as in Problem \ref{prob:cycle_selection} with 
  \begin{equation}
    (i\in\mathcal{I}) \ \uf{i}(\x) \coloneqq \frac{1}{2}\|x_i - x_{i-1}\|^2. 
  \end{equation}
  Case 2 satisfies Assumption~\ref{asm:G_monotone_Lipschitzian}\ref{item:asm_frakG} (see Proposition \ref{prop:uG_cycle}).
\end{enumerate}

\begin{figure}[t]
  \centering
  \begin{minipage}[b]{0.47\columnwidth}
      \centering
      \includegraphics[clip, width=0.9\columnwidth]{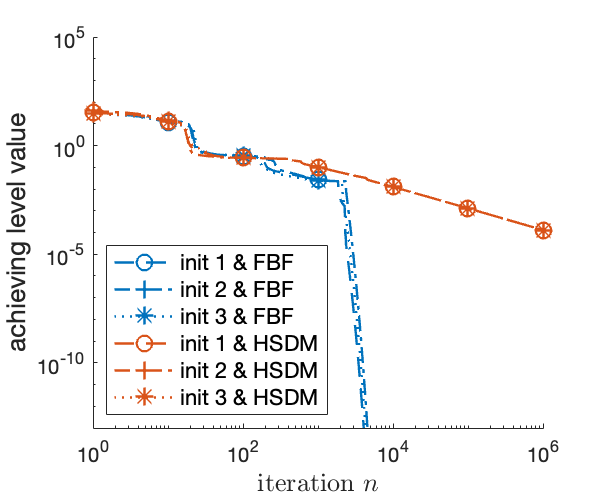}
      \label{fig:ex1.5_residual}
      \title{(Case 1)}
  \end{minipage}
  \hspace{0.01\columnwidth}
  \begin{minipage}[b]{0.485\columnwidth}
      \centering
      \includegraphics[clip, width=0.9\columnwidth]{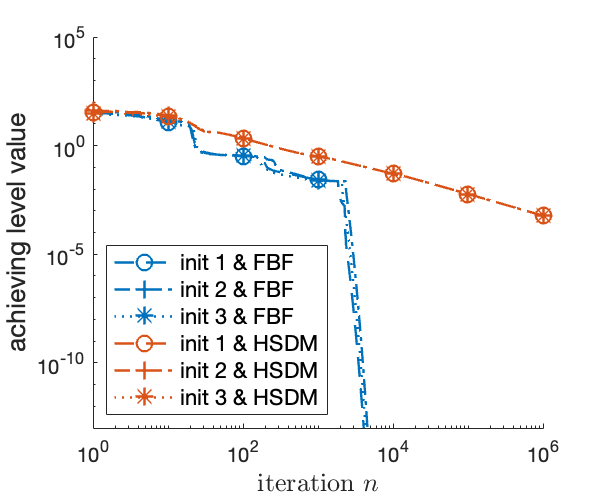}
      \label{fig:ex2_residual}
      \title{(Case 2)}
  \end{minipage}
  \caption{Achievement level $\|(P_{\overline{B}(0, r)} \circ \Talpha)(\bm{\xi}_n) - \bm{\xi}_n\|$ toward lower-level v-GNE at each iteration.}
  \label{fig:residual}
\end{figure}

We applied the algorithm \eqref{eq:replace_HSDM} (called ``HSDM'') with $(\gamma, \alpha, r) = (0.25, 0.75, 10^{15})$ to both cases. 
For comparison, we also used an iterative algorithm (called ``FBF''):
\begin{align}
  \label{eq:ex_FBF}
  (\forall n \in \N) \ \bm{\xi}_{n+1} = (P_{\overline{B}(0, r)} \circ \Talpha)(\bm{\xi}_n)
\end{align}
just for finding a solution of the lower-level v-GNEP because FBF can be seen as an instance of HSDM for $\uf{i} = 0 \ (i\in \mathcal{I})$ (i.e., FBF can be seen as a selection free version of HSDM). 
For both algorithms, we chose randomly 3 different initial points $\bxi_0$ (called ``init \{1, 2, 3\}''). 
To evaluate the achieving level at $\bxi_n$ toward the lower-level v-GNE, we used $\|(P_{\overline{B}(0, r)} \circ \Talpha)(\bm{\xi}_n) - \bm{\xi}_n\|$ motivated by the fact $\|(P_{\overline{B}(0, r)} \circ \Talpha)(\bm{\xi}_n) - \bm{\xi}_n\| = 0 \Leftrightarrow \bxi_n \in \Fix(P_{\overline{B}(0, r)} \circ \Talpha) = \overline{B}(0, r) \cap \Fix(\Talpha)$ (see also \eqref{eq:V_fixed}). 

\begin{figure}[t]
  \centering
  \includegraphics[clip, width=0.9\columnwidth]{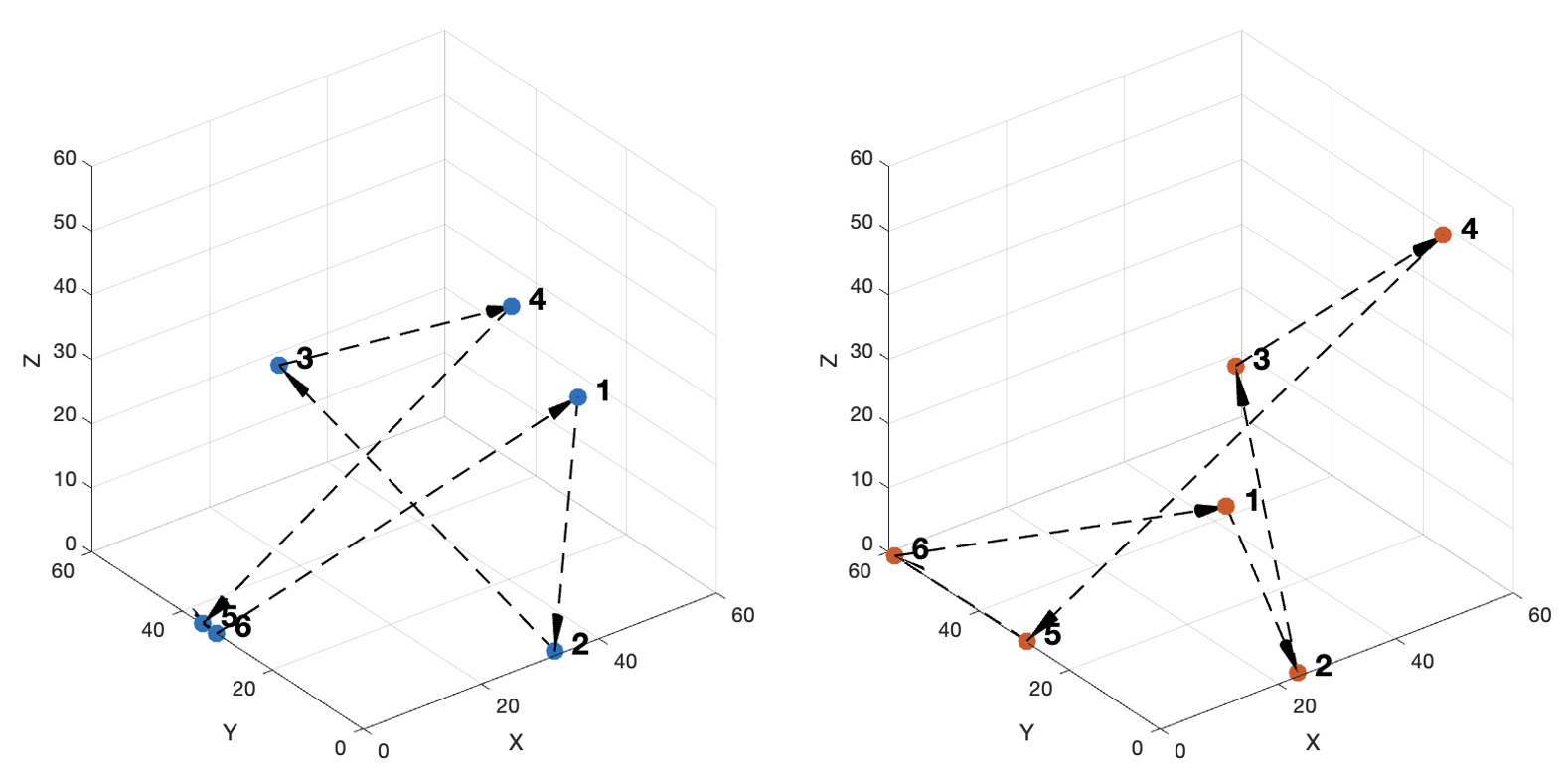}
  \caption{Visualization of points of all players' strategies in (Case 2) obtained by FBF (blue) as an approximation of a point in $\bm{\mathcal{V}}$, and obtained by HSDM (red) as an approximation of a cycle over $\bm{\mathcal{V}}$.}
  \label{fig:ex2_x}
\end{figure}

Fig. \ref{fig:residual} shows values of $\|(P_{\overline{B}(0, r)} \circ \Talpha)(\bm{\xi}_n) - \bm{\xi}_n\|$, in both (Case 1) and (Case 2), for each combination of initial point and algorithm. 
From Fig. \ref{fig:residual}, we can see that, in both cases, the value of $\|(P_{\overline{B}(0, r)} \circ \Talpha)(\bm{\xi}_n) - \bm{\xi}_n\|$ approaches zero, as $n$ increases, for every combination. 
\begin{figure}[t]
  \centering
  \includegraphics[clip, width=\columnwidth]{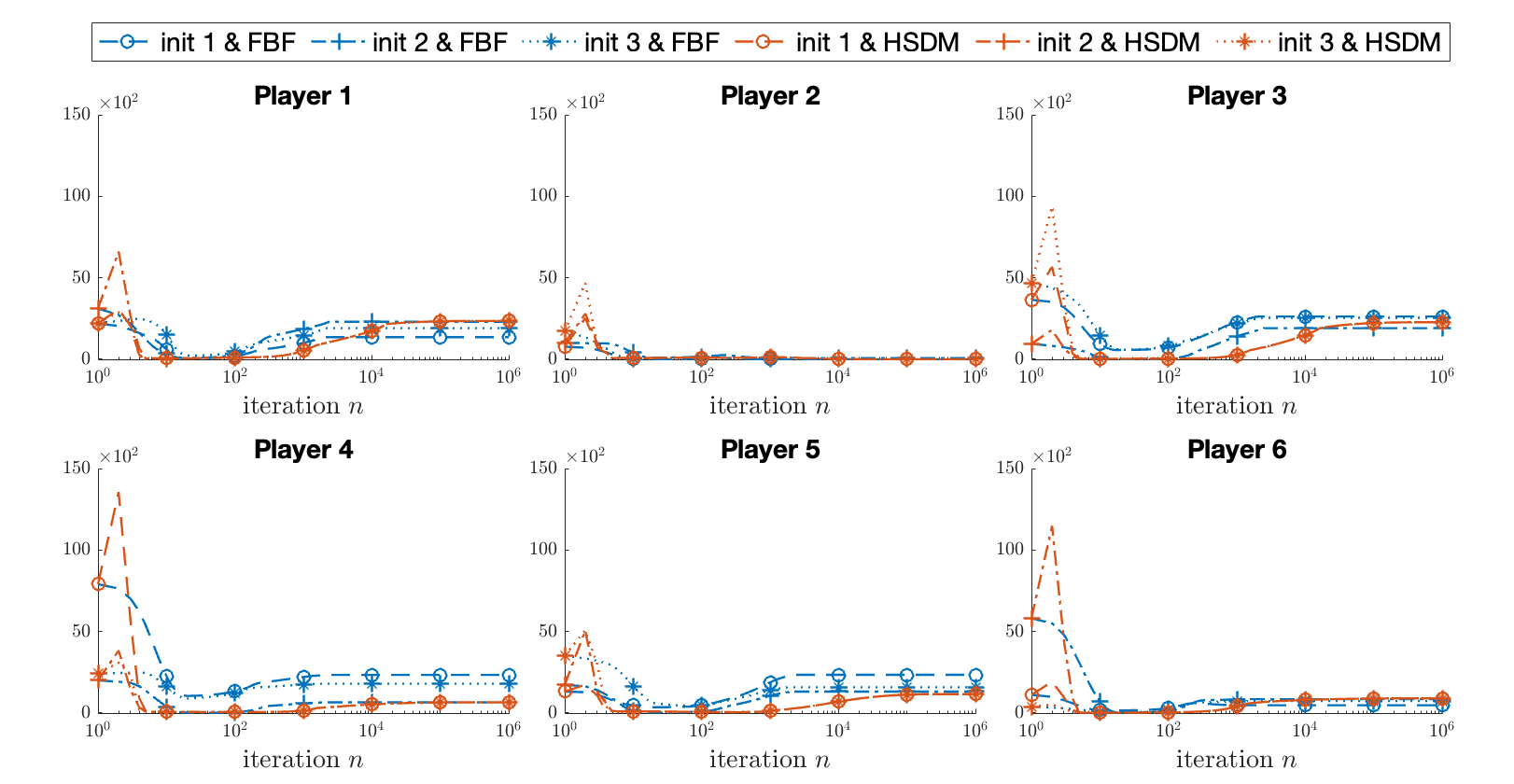}
  \caption{Each value of player $i(\in\mathcal{I})$'s upper level cost function (Case 1). }
  \label{fig:ex1.5_costUp}
\end{figure}
\begin{figure}[t]
  \centering
  \includegraphics[clip, width=\columnwidth]{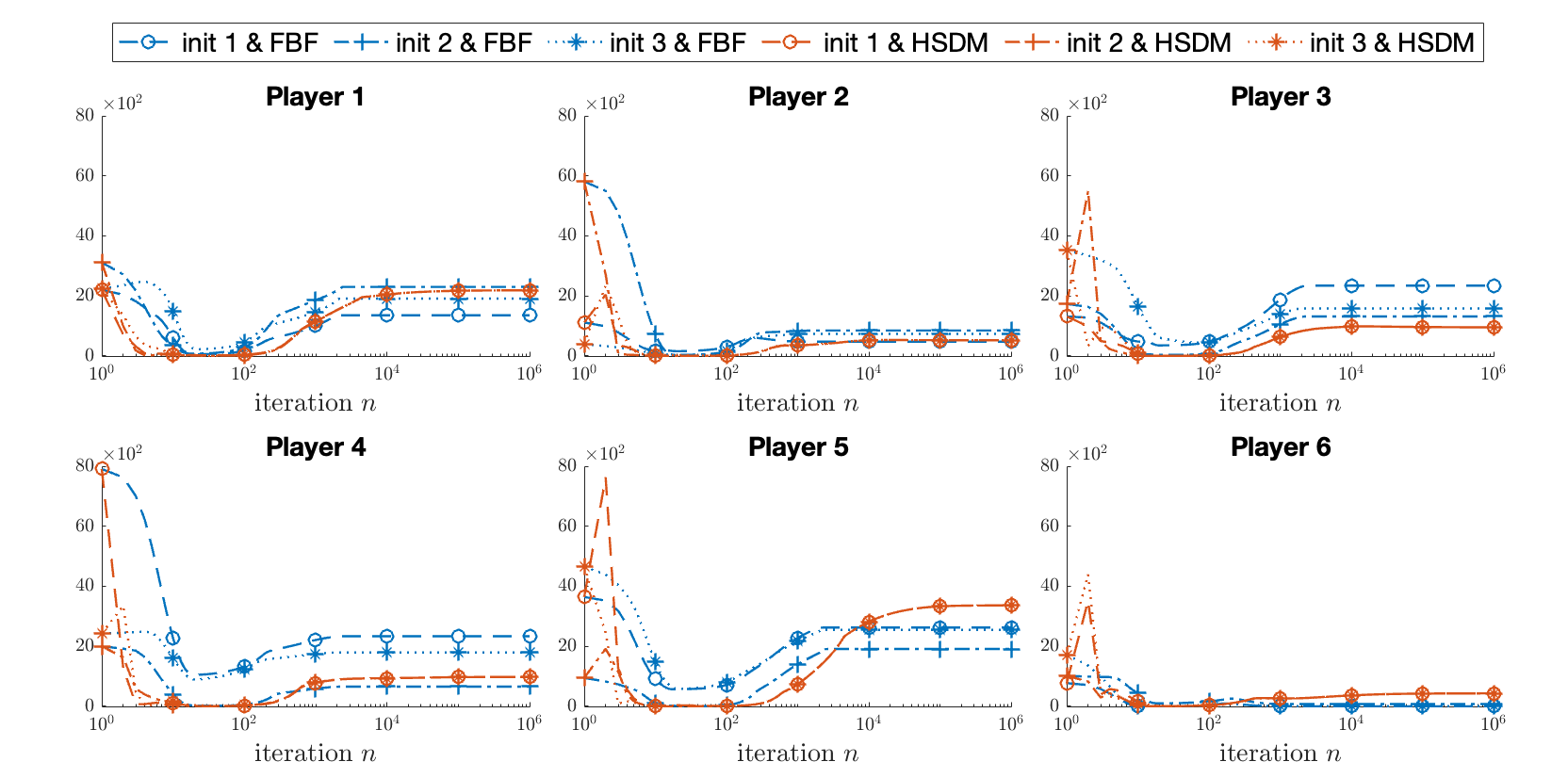}
  \caption{Each value of player $i(\in\mathcal{I})$'s upper level cost function (Case 2). }
  \label{fig:ex2_costUp}
\end{figure}
Next, we turn to examine that HSDM is making an effort for the further advanced goal than FBF (see Remark \ref{rem:HvGNEP}(ii)). 
Fig \ref{fig:ex1.5_costUp} and Fig. \ref{fig:ex2_costUp} show the values of upper level cost functions $\uf{i} \ (i\in\mathcal{I})$ in (Case 1) and (Case 2), respectively. 
From Fig \ref{fig:ex1.5_costUp} and Fig. \ref{fig:ex2_costUp}, we can see that: (i) $\x_n = \Q_{\bm{\mathcal{H}}}(\bxi_n)$ generated by FBF approaches different levels of $\uf{i} \ (i\in\mathcal{I})$ for different initial points, and (ii) $\x_n = \Q_{\bm{\mathcal{H}}}(\bxi_n)$ generated by HSDM approaches the same level of $\uf{i} \ (i\in\mathcal{I})$ for different initial points. 
This observation tells us that HSDM is achieving the advanced goal, i.e., equilibrium selection according to all players' upper-level requirements $\uf{i} \ (i\in\mathcal{I})$ while FBF achieves a lower-level v-GNE, depending on the choices of initial points, without any aimed selection. 
Fig. \ref{fig:ex2_x} visualizes the points of all players' strategies in (Case 2) obtained (blue) by FBF and (red) by HSDM, with a common initial point (init 2), after $10^6$ iterations. 
The points obtained by HSDM seem to achieve fairly balanced positions while the points obtained by FBF seem to have unexpected bias between $\|x_6  - x_5\|$ and other distances $\|x_i - x_{i-1}\|$. 

%% file: docs/Appendix.tex
\section{Known facts}
\label{sec:known_facts}
For readers' convenience, we present some known facts in convex analysis and monotone operator theory. 

\begin{fact}[Some properties of subdifferential]\label{fact:subdifferential}\,
	\begin{enumerate}[label=\textnormal{(\roman*)}]
		\item(Subdifferential and conjugate~\cite[Cor. 16.30]{BC2017}).
		      Let $f\in\Gamma_{0}(\mathcal{H})$.
		      Then, for any $(x,u)\in\mathcal{H}\times\mathcal{H}$, $u\in\partial f(x) \Leftrightarrow x\in \partial f^{\ast}(u)$, which implies that $(\partial f)^{-1} = \partial f^{\ast}$.
    \item(Subdifferential of separable sum~\cite[Prop. 16.9]{BC2017}).
          Let $(\mathcal{K}_1, \dots, \mathcal{K}_L)$ be finite-dimensional real Hilbert spaces, and let $f_l\in\Gamma_{0}(\mathcal{K}_l)$ for every $l \in \{1, \dots, L\}$. 
          Then $\partial \bigoplus_{l=1}^L f_l = \bigtimes_{l=1}^L \partial f_l$, where $\bigoplus_{l=1}^L f_l:\bigtimes_{l=1}^L\mathcal{K}_l \to (-\infty, \infty]:(z_1, \dots, z_L)\mapsto\sum_{l=1}^{L}f_l(z_l).$
    \item(Subdifferential of indicator function~\cite[Exm. 16.13]{BC2017}).
          For a nonempty closed convex set $C\subset\mathcal{H}$, $\partial \iota_C(x) = \{u\in\mathcal{H} \mid (\forall y\in C) \ \product[]{y-x}{u}\leq 0\}$ holds if $x\in C$. 
    \item \label{fact:chain_rule}
          (Chain rule~\cite[Cor. 16.53(i)]{BC2017}). 
          Let $g\in\Gamma_{0}(\mathcal{K})$ and $L\in\mathcal{B}(\mathcal{H}, \mathcal{K})$ satisfy $0_{\mathcal{K}} \in \operatorname{ri}\bigl(\dom(g) - \operatorname{ran}(L)\bigr)$.
          Then $\partial (g \circ L)=L^* \circ \partial g \circ L.$
  \end{enumerate}
\end{fact}

\begin{definition}[{Maximally monotone operator~\cite[Def. 20.20]{BC2017}}]
  \label{def:monotone_operator}
  A set-valued operator $\mathcal{A}:\mathcal{H}\to2^{\mathcal{H}}$ is said to be \emph{maximally monotone} if for every $(x, u)\in\mathcal{H}\times\mathcal{H}$, 
  \begin{equation}
    (x, u) \in \mathrm{gra}(\mathcal{A}) \Leftrightarrow 
    \bigl(\forall (y, v) \in \mathrm{gra}(\mathcal{A})\bigr) \ \product[\mathcal{H}]{x-y}{u-v} \geq 0,
  \end{equation}
  where $\mathrm{gra}(\mathcal{A})\coloneqq\{(x, u)\in\mathcal{H}\times\mathcal{H} \mid u\in\mathcal{A}(x)\}$. 
\end{definition}

\begin{fact}[Some properties of maximally monotone operators]
  \label{fact:maximally_monotone_operator} \ 
  \begin{enumerate}[label=\textnormal{(\roman*)}]
    \item (\cite[Prop. 20.23]{BC2017}). Let $\mathcal{K}_1$ and $\mathcal{K}_2$ be a finite-dimensional real Hilbert spaces, and let $\mathcal{A}_1:\mathcal{K}_1\to 2^{\mathcal{K}_1}$ and $\mathcal{A}_2:\mathcal{K}_2\to 2^{\mathcal{K}_2}$ be maximally monotone operators. 
    Then $\bm{\mathcal{A}}:\mathcal{K}_1\times\mathcal{K}_2\to 2^{\mathcal{K}_1\times\mathcal{K}_2}:(x_1, x_2)\mapsto\mathcal{A}_1(x_1)\times\mathcal{A}_2(x_2)$ is maximally monotone. 
    \item (\cite[Thm. 20.25]{BC2017}). The subdifferential $\partial f:\mathcal{H}\to 2^{\mathcal{H}}$ of $f\in\Gamma_0(\mathcal{H})$ is maximally monotone. 
    \item \label{item:single_resolvent} (\cite[Cor. 23.11(i)]{BC2017}) Let $\mathcal{A}:\mathcal{H}\to2^{\mathcal{H}}$ be a maximally monotone operator. 
    For any $\gamma > 0$, $(\Id + \gamma \mathcal{A})^{-1}$ is $\frac{1}{2}$-averaged nonexpansive\footnote{This implies that $(\Id + \gamma \mathcal{A})^{-1}$ is single-valued.}, i.e., for every $x \in \mathcal{H}$ and $y \in \mathcal{H}$, $\|(\Id + \gamma \mathcal{A})^{-1}(x) - (\Id + \gamma \mathcal{A})^{-1}(y)\|^2_{\mathcal{H}} \leq \|x-y\|^2_{\mathcal{H}} - \bigl\|(x - y) - \bigl( (\Id + \gamma \mathcal{A})^{-1}(x) - (\Id + \gamma \mathcal{A})^{-1}(y)\bigr)\bigr\|^2_{\mathcal{H}}$. 
  \end{enumerate}
\end{fact}

\section{Proof of Proposition \ref{prop:FBF_fix}}
\label{Appendix:proof_prop_FBF_fix}

\underline{Proof of (i)}: Assumption \ref{asm:G_monotone_Lipschitzian}\ref{item:asm_sum_chain} yields $\partial \iota_{\bm{\mathfrak{C}}} = \partial (\iota_{\C} + \iota_D \circ \Li) = \partial \iota_{\C} + \Li^* \circ \partial \iota_D \circ \Li$.  
By noting that $\iota_{\C}$ is the separable sum of $\iota_{C_i} \ (i \in \mathcal{I})$, Fact \ref{fact:subdifferential}(ii) gives $\partial \iota_{\C} = \bigtimes_{i\in \mathcal{I}} \partial \iota_{C_i}$. 
Since $\bm{\mathcal{V}}$ in \eqref{eq:VE} is the solution set of the variational inequality $\VI(\bm{\mathfrak{C}}, \G)$, 
\begin{equation}
  \label{eq:VI_normalcone}
  \begin{split}
    \x \in \bm{\mathcal{V}} = \left\{\vv \in \bm{\mathfrak{C}} \mid (\forall \w \in \bm{\mathfrak{C}}) \ \product[\bm{\mathcal{H}}]{\G(\vv)}{\w-\vv} \geq 0 \right\} 
    \Leftrightarrow -\G(\x) \in  \partial \iota_{\bm{\mathfrak{C}}}(\x)
    \Leftrightarrow \bm{0}_{\bm{\mathcal{H}}} \in  \partial \iota_{\bm{\mathfrak{C}}}(\x) + \G(\x)
  \end{split}
\end{equation}
holds true, where the first relation follows from Fact \ref{fact:subdifferential}(iii). 
By substituting $\partial \iota_{\bm{\mathfrak{C}}} = \partial \iota_{\C} + \Li^* \circ \partial \iota_D \circ \Li$ and $\partial \iota_{\C} = \bigtimes_{i\in \mathcal{I}} \partial \iota_{C_i}$, we have 
\begin{align}
  \x \in \bm{\mathcal{V}}
  &\Leftrightarrow \bm{0}_{\bm{\mathcal{H}}} \in \bigtimes_{i\in I} \partial \iota_{C_i}(x_i) + \Li^* \partial \iota_D (\Li\x) + \G(\x) \\
  &\Leftrightarrow (\exists u \in \mathcal{G}) \ 
  \begin{cases}
    \bm{0}_{\bm{\mathcal{H}}} \in \displaystyle\bigtimes_{i\in I} \partial \iota_{C_i}(x_i) + \Li^* u + \G(\x) \\
    u \in \partial \iota_D (\Li\x)
  \end{cases} \\
  &\hspace{-1.6em}\overset{\mathrm{Fact} \ \mathrm{\ref{fact:subdifferential}(i)}}{\Leftrightarrow} (\exists u \in \mathcal{G}) \ 
  \begin{cases}
    \bm{0}_{\bm{\mathcal{H}}} \in \displaystyle\bigtimes_{i\in I} \partial \iota_{C_i}(x_i) + \Li^* u + \G(\x) \\
    \Li\x \in \partial \iota_D^* (u)
  \end{cases} \\
  \label{eq:VI_zero_1}
  &\Leftrightarrow (\exists u \in \mathcal{G}) \ 
  \begin{cases}
    \bm{0}_{\bm{\mathcal{H}}} \in \displaystyle\bigtimes_{i\in I} \partial \iota_{C_i}(x_i) + \Li^* u + \G(\x) \\
    \bm{0}_{\mathcal{G}} \in \partial \iota_D^* (u) - \Li\x
  \end{cases} 
  \Leftrightarrow 
  (\exists u \in \mathcal{G}) \ 
  \bm{0}_{\bm{\mathcal{H}}\times\mathcal{G}} \in \A(\x, u) + \B(\x, u). 
\end{align}
By the assumption $\bm{\mathcal{V}}\neq \varnothing$ (see Problem \ref{prob:H_vGNEP}), we have $\zer(\A+\B) \neq \varnothing$, and thus $\bm{\mathcal{V}} = \Q_{\bm{\mathcal{H}}}\bigl( \zer(\A+\B) \bigr)$. 

As the final step of proof of (i), $\zer(\A+\B) = \Fix(\TFBF)$ can be verified as follows (a similar discussion is found in \cite{Franci2020FBFforGeneralizedNash}). 

(Proof of $\zer(\A+\B) \subset \Fix(\TFBF)$) 
Let $\bm{\xi} \in \zer(\A+\B)$. 
Since $\B$ is maximally monotone (see Definition \ref{def:monotone_operator} and Lemma \ref{lemma:A_B_monotone}), $\zer(\A+\B)$ can be expressed as \cite[Proposition 26.1(iv)(a)]{BC2017}, 
\begin{equation}
  \label{eq:zer_fix_FB}
  \zer(\A+\B) = \Fix(\T_{\FB}), 
\end{equation}
where the forward-backward operator $\T_{\FB}:\bm{\mathcal{H}}\times\mathcal{G}\to\bm{\mathcal{H}}\times\mathcal{G}$ is defined in \eqref{eq:FBF}. 
A straightforward manipulation from \eqref{eq:zer_fix_FB} yields
\begin{align}
  &\bm{\xi} \in \zer(\A+\B) 
  \Leftrightarrow \bxi = \T_{\FB}(\bxi)
  \Rightarrow \bxi - \gamma \A(\bxi) = \T_{\FB}(\bxi) - \gamma \A \bigl( \T_{\FB}(\bxi) \bigr) \\
  &\Leftrightarrow \bxi = (\Id - \gamma \A) \bigl( \T_{\FB}(\bxi) \bigr) + \gamma \A(\bxi) = \TFBF(\bxi)
  \Leftrightarrow \bxi \in \Fix(\TFBF). 
\end{align}

(Proof of $\zer(\A+\B) \supset \Fix(\TFBF)$) 
Assume contrarily $\zer(\A+\B) \not\supset \Fix(\TFBF)$. 
By $\zer(\A+\B) \subset \Fix(\TFBF)$, there exists some $\bxi \in \Fix(\TFBF) \setminus \zer(\A+\B)$. 
Since $\TFBF$ in \eqref{eq:FBF} can be expressed as $\TFBF = (\Id - \gamma\A)\circ\T_{\FB} + \gamma\A$, we have $\bxi - \T_{\FB}(\bxi)\bigl(=\TFBF(\bxi) - \T_{\FB}(\bxi)\bigr) = \gamma \bigl( \A(\bxi) - \A \bigl( \T_{\FB}(\bxi) \bigr) \bigr)$. 
By noting that $\A$ is $(\kappa_{\G} + \|\Li\|_{\mathrm{op}})$-Lipschitzian over $\bm{\mathcal{H}}\times\mathcal{G}$ (see Lemma \ref{lemma:A_B_monotone}), we have 
\begin{align}
  \label{eq:xi_FB}
  \|\bxi - \T_{\FB}(\bxi)\| 
  = \gamma \Bigl\|\A(\bxi) - \A \bigl( \T_{\FB}(\bxi) \bigr)\Bigr\| 
  \leq \gamma (\kappa_{\G} + \|\Li\|_{\mathrm{op}}) \|\bxi - \T_{\FB}(\bxi)\|. 
\end{align}
Since $\gamma \in \bigl(0, 1/(\kappa_{\G} + \|\Li\|_{\mathrm{op}})\bigr)$ and $\bxi \notin \zer(\A+\B) = \Fix(\T_{\mathrm{FB}}) \ (\mathrm{i.e.,} \ \|\bxi - \T_{\mathrm{FB}}(\bxi)\| \neq 0)$, we obtain $\|\bxi - \T_{\FB}(\bxi)\| < \|\bxi - \T_{\FB}(\bxi)\|$ from \eqref{eq:xi_FB}, which is absurd.

Therefore, we have $\Fix(\T_{\FBF}) = \zer(\A+\B)$. 
Moreover, by $\Fix(\T_{\FBF}^{\alpha} )= \Fix(\TFBF)$ (see Definition \ref{df:quasi_nonexpansive}(v)), we complete the proof of (i). 

\noindent\underline{Proof of (ii)}: Firstly, we show that $\TFBF$ is quasi-nonexpansive. 
Recall that $\A$ and $\B$ are maximally monotone (see Lemma \ref{lemma:A_B_monotone}) and $\zer(\A+\B) = \Fix(\TFBF) \neq \varnothing$ (see Proposition \ref{prop:FBF_fix}(i)). 
By defining $\T_{\FB}$ as in \eqref{eq:FBF}, \cite[Lemma 3.1]{tsengModifiedForwardBackwardSplitting2000} yields 
\begin{align}
  &(\forall \gamma \in (0, 1/(\kappa_{\G} + \|\Li\|_{\mathrm{op}})), \forall \bxi \in \dom(\A)=\bm{\mathcal{H}}\times\mathcal{G}, \forall \bzeta \in \zer(\A+\B) = \Fix(\TFBF), \exists \eta \geq 0) \\
  & \qquad \qquad \|\TFBF(\bxi) - \bzeta\|^2 = \|\bxi - \bzeta\|^2 + \gamma^2 \Bigl\|\A\bigl(\T_{\FB}(\bxi)\bigr) - \A(\bxi)\Bigr\|^2 - \|\T_{\FB}(\bxi) - \bxi\|^2 - 2\gamma \eta. 
\end{align}
Since $\A$ is $(\kappa_{\G} + \|\Li\|_{\mathrm{op}}=:\kappa_{\A})$-Lipschitzian over $\bm{\mathcal{H}}\times\mathcal{G}$, we get an upper bound 
\begin{equation}
  \label{eq:FBF_quasi}
  \begin{aligned}
    &\bigl( \forall \bxi \in \bm{\mathcal{H}}\times\mathcal{G}, \forall \bzeta \in \Fix(\TFBF) \bigr) \quad \\
    &\qquad \qquad  \|\TFBF(\bxi) - \bzeta\|^2 \leq \|\bxi - \bzeta\|^2 - (1-\gamma^2 \kappa_{\A}^2) \|\T_{\FB}(\bxi) - \bxi\|^2 \leq \|\bxi - \bzeta\|^2, 
  \end{aligned}
\end{equation}
where the last inequality follows from $\gamma \in (0,\kappa_{\A}^{-1})$. 
Then, $\TFBF$ is quasi-nonexpansive, and thus by Fact \ref{fact:quasi_nonexpansive}(ii), the $\alpha$-averaged operator $\T_{\FBF}^{\alpha}$ in \eqref{eq:Talpha} is strongly attracting. 

Moreover, since $(\Id + \gamma \B)^{-1}$ is continuous from the maximal monotonicity of $\B$ (see Fact \ref{fact:maximally_monotone_operator}\ref{item:single_resolvent}), $\TFBF$ is continuous, and thus $\T_{\FBF}^{\alpha}$ is also continuous. 
The continuity of $\T_{\FBF}^{\alpha}$ implies that $\T_{\FBF}^{\alpha} - \Id$ is demi-closed at $\bm{0}_{\bm{\mathcal{H}}\times\mathcal{G}}$. 
Therefore, Fact \ref{fact:quasi_shrinking_stquasi} guarantees that $\T_{\FBF}^{\alpha}$ is quasi-shrinking on any bounded closed convex set $\bm{\mathcal{C}}\subset\bm{\mathcal{H}}\times\mathcal{G}$ satisfying $\Fix(\T_{\FBF}^{\alpha})\cap\bm{\mathcal{C}} \neq \varnothing$. 
\qed

\section{Properties of $\A$ in \eqref{eq:A} and $\B$ in \eqref{eq:B}}
\begin{lemma}
  \label{lemma:A_B_monotone}
  Under Assumption \ref{asm:G_monotone_Lipschitzian}, 
  $\A$ in \eqref{eq:A} is monotone and $\kappa_{\A}:=(\kappa_{\G} + \|\Li\|_{\mathrm{op}})$-Lipschitzian over $\bm{\mathcal{H}}\times\mathcal{G}$ (which imply that $\A$ is maximally monotone \cite[Corollary 20.28]{BC2017}), and $\B$ in \eqref{eq:B} is maximally monotone (see Definition \ref{def:monotone_operator} for maximal monotonicity). 
\end{lemma}
\textit{Proof.} By the monotonicity of $\G$, we have for any $(\x, u), (\y, v) \in \bm{\mathcal{H}} \times \mathcal{G}$, 
\begin{align}
  & \product[\bm{\mathcal{H}} \times \mathcal{G}]{(\x, u) - (\y, v)}{\A(\x, u) - \A(\y, v)} \\
  &= \product[\bm{\mathcal{H}}]{\x - \y}{\G(\x) - \G(\y)} + \product[\bm{\mathcal{H}}]{\x -\y}{\Li^*u - \Li^*v} 
  - \product[\mathcal{G}]{u - v}{\Li \x - \Li \y}\\
  &= \product[\bm{\mathcal{H}}]{\x - \y}{\G(\x) - \G(\y)} + \product[\mathcal{G}]{\Li \x - \Li \y}{u - v}
  - \product[\mathcal{G}]{u - v}{\Li \x - \Li \y}\\
  &= \product[\bm{\mathcal{H}}]{\x - \y}{\G(\x) - \G(\y)} 
  \geq 0,  
\end{align}
which implies the monotonicity of $\A$. 
To show that $\A$ is Lipschitzian, consider the decomposition $\A=\A_1 + \A_2$ with $\A_1 : \bm{\mathcal{H}} \times\mathcal{G} \rightarrow \bm{\mathcal{H}} \times\mathcal{G}:
(\x, u)\mapsto(\G(\x), 0_{\mathcal{G}})$ 
and $\A_2 : \bm{\mathcal{H}} \times\mathcal{G} \rightarrow \bm{\mathcal{H}} \times\mathcal{G}:(\x, u)
\mapsto(\Li^*u, -\Li\x)$. Since $\G$ is $\kappa_{\G}$-Lipschitzian over $\bm{\mathcal{H}}$, for every $(\x, u), (\y, v) \in \bm{\mathcal{H}} \times \bm{\mathcal{G}}$, we have $\|\A_1(\x, u) - \A_1(\y, v)\|_{\bm{\mathcal{H}} \times \mathcal{G}} = \|\G(\x) - \G(\y)\|_{\bm{\mathcal{H}}} \leq \kappa_{\G}\|\x - \y\|_{\bm{\mathcal{H}}}$. 
We also have $\|\A_2(\x, u) - \A_2(\y, v)\|_{\bm{\mathcal{H}} \times \mathcal{G}} = \sqrt{\|\Li^*u - \Li^*v\|_{\bm{\mathcal{H}}}^2 + \|\Li\x - \Li\y\|_{\mathcal{G}}^2} \leq \|\Li\|_{\mathrm{op}} \sqrt{\|\x -\y\|_{\bm{\mathcal{H}}}^2 + \|u - v\|_{\mathcal{G}}^2} 
= \|\Li\|_{\mathrm{op}} \|(\x, u) - (\y, v)\|_{\bm{\mathcal{H}} \times \mathcal{G}}$. 
From these inequalities, for every $(\x, u), (\y, v) \in \bm{\mathcal{H}} \times \bm{\mathcal{G}}$, we get 
\begin{align}
  \|\A(\x, u) - \A(\y, v)\|_{\bm{\mathcal{H}} \times \mathcal{G}}
  & = \|\A_1(\x, u) - \A_1(\y, v) + \A_2(\x, u) - \A_2(\y, v)\|_{\bm{\mathcal{H}} \times \mathcal{G}} \\
  & \leq \|\A_1(\x, u) - \A_1(\y, v)\|_{\bm{\mathcal{H}} \times \mathcal{G}} + \|\A_2(\x, u) - \A_2(\y, v)\|_{\bm{\mathcal{H}} \times \mathcal{G}} \\
  &\leq \kappa_{\G}\|\x - \y\|_{\bm{\mathcal{H}}} + 
  \|\Li\|_{\mathrm{op}} \|(\x, u) - (\y, v)\|_{\bm{\mathcal{H}} \times \mathcal{G}} \\
  & \leq (\kappa_{\G} + \|\Li\|_{\mathrm{op}}) \|(\x, u) - (\y, v)\|_{\bm{\mathcal{H}} \times \mathcal{G}}, 
\end{align}
which implies that $\A$ is $(\kappa_{\G} + \|\Li\|_{\mathrm{op}})$-Lipschitz continuous over $\bm{\mathcal{H}}\times\mathcal{G}$. 

From $\iota_{C_i} \in \Gamma_0(\mathcal{H}_i) \ (\forall i \in \mathcal{I})$, and $\iota_D^{\ast} \in \Gamma_0(\mathcal{G})$ by \cite[Cor. 13.38]{BC2017}, $\partial \iota_{C_i} \ (i \in \mathcal{I})$ and $\partial \iota_D^{\ast}$ are maximally monotone by Fact \ref{fact:maximally_monotone_operator}(ii). 
Hence, the operator $\B$ is maximally monotone by Fact \ref{fact:maximally_monotone_operator}(i). 
\qed

\section{Proof of Theorem \ref{thm:HSDM_convergence}}
To invoke Fact \ref{fact:HSDM}, we check below that (a) $\T_{\FBF}^{\alpha}$ is a quasi-nonexpansive operator with bounded $\Fix(\T_{\FBF}^{\alpha})$; (b) $\widetilde{\uG}$ is paramonotone over $\Fix(\Talpha)$ and Lipschitzian over $\Talpha(\bm{\mathcal{H}}\times\mathcal{G})$; (c) there exists some nonempty bounded closed convex set $\bm{K} \subset \bm{\mathcal{H}}\times\mathcal{G}$ such that $(\bxi_n)_{n\in\N} \subset \bm{K}$ and $\T_{\FBF}^{\alpha}$ is quasi-shrinking on $\bm{K}$. 

(a) See Proposition \ref{prop:FBF_fix}(ii) and Condition \eqref{eq:bound_cond_proto}. 

(b) By using $\Q_{\bm{\mathcal{H}}}^*:\bm{\mathcal{H}}\to\bm{\mathcal{H}}\times\mathcal{G}:\x\mapsto(\x, \bm{0}_{\mathcal{G}})$ for the canonical projection $\Q_{\bm{\mathcal{H}}}$ onto $\bm{\mathcal{H}}$ in Proposition \ref{prop:FBF_fix}(i), $\widetilde{\uG}$ can be expressed as $\widetilde{\uG} = \Q_{\bm{\mathcal{H}}}^*\circ\uG\circ\Q_{\bm{\mathcal{H}}}$. 
Then, \cite[Proposition 22.2(ii)]{BC2017} guarantees that $\widetilde{\uG}$ is paramonotone over $\bm{\mathcal{H}}\times\mathcal{G}$. 
In addition, since $\uG$ is $\kappa_{\bm{\mathfrak{G}}}$-Lipschitzian over $\bm{\mathcal{H}}$, we have for every $(\x, u)\in\bm{\mathcal{H}}\times\mathcal{G}$ and every $(\y, v)\in\bm{\mathcal{H}}\times\mathcal{G}$, 
\begin{align}
  \|\widetilde{\uG}(\x, u) - \widetilde{\uG}(\y, v)\|^2_{\bm{\mathcal{H}}\times\mathcal{G}}
  = \|\uG(\x) - \uG(\y)\|_{\bm{\mathcal{H}}}^2 
  \leq \kappa_{\bm{\mathfrak{G}}}^2 \|\x - \y\|_{\bm{\mathcal{H}}}^2 
  \leq \kappa_{\bm{\mathfrak{G}}}^2 \|\left(\x, u\right) - \left(\y, v\right)\|^2_{\bm{\mathcal{H}}\times\mathcal{G}}, 
\end{align}
which implies that $\widetilde{\uG}$ is Lipschitzian over $\bm{\mathcal{H}}\times\mathcal{G}$. 

(c) By Condition \eqref{eq:bound_cond_proto}, there exists a closed ball $\clBall\coloneqq\{\bxi\in\bm{\mathcal{H}}\mid \|\bxi\|_{\bm{\mathcal{H}}
\times\mathcal{G}} \leq r\} \subset \bm{\mathcal{H}}\times\mathcal{G}$ satisfying $(\bxi_n)_{n\in\N} \subset \clBall$ and $\Fix(\T_{\FBF}^{\alpha}) \subset \clBall$. 
Then, Proposition \ref{prop:FBF_fix}(ii) guarantees that $\T_{\FBF}^{\alpha}$ is quasi-shrinking on $\clBall$. 

By applying Fact \ref{fact:HSDM}, we get $\lim_{n \to \infty} d\bigl(\bxi_n, \bm{\Omega} \bigr) = 0$. 
Since $\Q_{\bm{\mathcal{H}}}\bigl(P_{\bm{\Omega} }(\bxi)\bigr) \in \uV$ holds for any $\bxi \in \bm{\mathcal{H}}\times\mathcal{G}$, we deduce that 
\begin{equation}
  \label{eq:dist_spaces}
  \begin{split}
    d(\x_n, \uV) 
    &=\Bigl\|\x_n - P_{\uV}(\x_n)\Bigr\|_{\bm{\mathcal{H}}} 
    \leq \biggl\|\x_n - \Q_{\bm{\mathcal{H}}}\left(P_{\bm{\Omega}}(\bxi_n)\right)\biggr\|_{\bm{\mathcal{H}}} \\
    &= \biggl\|\Q_{\bm{\mathcal{H}}}(\bxi_n) - \Q_{\bm{\mathcal{H}}}\left(P_{\bm{\Omega}}(\bxi_n)\right)\biggr\|_{\bm{\mathcal{H}}} 
    \leq \Bigl\|\bxi_n - P_{\bm{\Omega} }(\bxi_n)\Bigr\|_{\bm{\mathcal{H}}\times\mathcal{G}}
    = d\bigl(\bxi_n, \bm{\Omega} \bigr).
  \end{split}
\end{equation}
We obtain $\lim_{n \rightarrow \infty} d(\x_n, \uV) = 0$ by $\lim_{n \rightarrow \infty} d\bigl(\bxi_n, \bm{\Omega}\bigr) = 0$. 

\section{Proof of Theorem \ref{thm:modified_HSDM}}
\underline{Proof of (i)}: Since, $P_{\clBall}$ is a $\frac{1}{2}$-averaged nonexpansive operator \cite[Cor. 4.18]{BC2017}, $P_{\clBall}$ is continuous and $\frac{1}{2}$-averaged quasi-nonexpansive with bounded $\Fix\bigl(P_{\clBall}\bigr) = \clBall \neq \varnothing$. 
Then, $P_{\clBall}$ is strongly attracting (see Fact \ref{fact:quasi_nonexpansive}(ii)) and continuous. 
Since $\Talpha$ is continuous and strongly attracting (see Proposition \ref{prop:FBF_fix}(ii)), $P_{\clBall}\circ\Talpha$ is continuous and strongly attracting with $\Fix\bigl(P_{\clBall}\circ\Talpha\bigr) = \Fix\bigl( \Talpha \bigr)\cap\Fix\bigl(P_{\clBall}\bigr) = \Fix\bigl( \Talpha \bigr)\cap\clBall = \Fix\bigl( \Talpha \bigr)$ (see Fact \ref{fact:quasi_nonexpansive}(iii)).
Since the continuity of $P_{\clBall}\circ\Talpha$ implies the demi-closedness at $\bm{0}_{\bm{\mathcal{H}}\times\mathcal{G}}$ of $P_{\clBall}\circ\Talpha$, $P_{\clBall}\circ\Talpha$ is quasi-shrinking on any bounded closed convex set $\bm{\mathcal{C}}\subset\bm{\mathcal{H}}\times\mathcal{G}$ satisfying $\Fix\bigl( P_{\clBall}\circ\Talpha \bigr)\cap\bm{\mathcal{C}} \neq \varnothing$ (see Fact \ref{fact:quasi_shrinking_stquasi}). 

\noindent\underline{Proof of (ii)}: The continuity of $\widetilde{\uG}$ and compactness of $\clBall$ imply \linebreak $\sup_{\bxi \in \clBall} \|\widetilde{\uG}(\bxi)\| \leq \tau_1$ with some $\tau_1 > 0$. 
Since we have $\bigl( \bigl( P_{\clBall}\circ\Talpha \bigr)(\bxi_n) \bigr)_{n\in\N} \subset \clBall$ and $\sup_{n \in \N}|\lambda_n| \leq \tau_2$ with some $\tau_2 > 0$ by (H1) in \eqref{eq:H1_H2}, $(\bxi_n)_{n\in\N}$ is bounded from $\|\bxi_{n+1}\| \leq \bigl\|\bigl( P_{\clBall}\circ\Talpha \bigr)(\bxi_n)\bigr\|+|\lambda_{n+1}|\biggl\|\widetilde{\uG}\Bigl(\bigl(P_{\clBall}\circ\Talpha\bigr)(\bxi_n)\Bigr)\biggr\| \leq r + \tau_1 \tau_2$ for all $n\in\N$. 

\noindent\underline{Proof of (iii)}: 
By $\Fix\bigl(P_{\clBall}\circ\Talpha\bigr) = \Fix\bigl( \Talpha \bigr) \cap \clBall = \Fix(\Talpha)$, $\bm{\Omega}$ in \eqref{eq:upper_VI_over_fixed_point} can be expressed as 
\begin{equation}
  \label{eq:replace_upper_VI_over_fixed_point_2}
  \hspace*{-0.8em}
  \begin{aligned}
    \bm{\Omega}= \Bigl\{ \bm{\xi} \in \Fix\bigl(P_{\clBall}\circ\Talpha\bigr) \Bigm\vert \bigl(\forall \bm{\zeta} \in \Fix\bigl(P_{\clBall}\circ\Talpha\bigr) \bigr) \ \bigl\langle \widetilde{\bm{\mathfrak{G}}}(\bm{\xi}), \bm{\zeta} - \bm{\xi} \bigr\rangle_{\bm{\mathcal{H}}\times\mathcal{G}} \geq 0 \Bigr\}. 
  \end{aligned}
\end{equation} 
Recall that (a) $P_{\clBall}\circ\Talpha$ is a quasi-nonexpansive operator with bounded $\Fix(P_{\clBall}\circ\Talpha)$ (see (i)), (b) $\widetilde{\uG}$ is paramonotone and Lipschitzian over $\bm{\mathcal{H}}\times\mathcal{G}$ (see the proof of Theorem \ref{thm:HSDM_convergence}(b) in Appendix D), and (c) there exists some nonempty bounded closed convex set $\bm{K} \subset \bm{\mathcal{H}}\times\mathcal{G}$ such that $(\bxi_n)_{n\in\N} \subset \bm{K}$ (see (ii)) and $P_{\clBall}\circ\Talpha$ is quasi-shrinking on $\bm{K}$ (see (i)). 
Then, we have $\lim_{n \rightarrow \infty} d\bigl(\bxi_n, \bm{\Omega} \bigr) = 0$ by applying Fact \ref{fact:HSDM} to \eqref{eq:replace_upper_VI_over_fixed_point_2} and $(\bxi_n)_{n\in\N} = (\x_n, u_n)_{n\in\N}$. 
Similarly to \eqref{eq:dist_spaces} in Appendix D, we have $\lim_{n \rightarrow \infty} d\bigl(\x_n, \uV \bigr) = 0$. 
\qed

\section{Proof of Proposition \ref{prop:bound_cond}}
Assume contrarily $\Fix(\Talpha)$ is unbounded. 
Then, there exists a sequence $(\x_n, u_n)_{n\in\N} \subset \Fix(\Talpha)$ such that $\|(\x_n, u_n)\| \rightarrow \infty \ (n\to \infty)$. 
We have $M \coloneqq \sup_{\x\in\bm{\mathcal{V}}}(\|\x\|) < \infty$ by \ref{item:bound_cond2}. 
From $\x_n \in \bm{\mathcal{V}}$ by \eqref{eq:V_fixed}, we have $\|(\x_n, u_n)\|^2 = \|\x_n\|^2 + \|u_n\|^2 \leq M^2 + \|u_n\|^2$, which implies $\|u_n\| \rightarrow \infty \ (n \to \infty)$. 
From \ref{item:bound_cond3}, we have $\|\Li^*u_n\| \to \infty \ (n\to\infty)$ \footnote{
Let $\lambda_{\min}$ be the minimum eigenvalue of $\Li\Li^*$. 
Since $\Li\Li^*$ is positive definite by \ref{item:bound_cond3}, for every $u\in\mathcal{G}$ such that $\|u\| \neq 0$, we have $0 < \lambda_{\min} \leq \frac{\product[]{\Li\Li^*u}{u}}{\product[]{u}{u}} = \frac{\|\Li^*u\|^2}{\|u\|^2}$, and consequently $0 < \lambda_{\min}\|u\|^2 \leq \|\Li^*u\|^2$. 
By $\|u_n\| \rightarrow \infty$, we have $\|\Li^*u_n\| \to \infty \ (n\to\infty)$.}. 
In addition, we have 
\begin{equation}
  \label{eq:-G-L}
  \|-\G(\x_n) - \Li^*u_n\| \geq \|\Li^*u_n\| - \|\G(\x_n)\| \geq \|\Li^*u_n\| - M_{\G}, 
\end{equation}
where $M_{\G} \coloneqq \sup_{\x \in \bm{\mathcal{V}}} (\|\G(\x)\|) < \infty$ by the compactness of $\bm{\mathcal{V}}$ and the continuity of $\G$. 
By $\|\Li^*u_n\| \to \infty \ (n\to\infty)$, we have $\|-\G(\x_n) - \Li^*u_n\| \to \infty \ (n\to\infty)$. 
Hence, there exists some subsequence $(\x_{n(k)}, u_{n(k)})_{k\in\N}$ of $(\x_n, u_n)_{n\in\N}$ such that $\|\Li^* u_{n(k)}\| \neq 0$ and $\|-\G(\x_{n(k)}) - \Li^* u_{n(k)} \| \neq 0$. 
Define $(\vv_k)_{k\in\N} \subset \bm{\mathcal{H}}$~and~$(\w_k)_{k\in\N} \subset \bm{\mathcal{H}}$~by
\begin{equation}
  \vv_k = \frac{-\G(\x_{n(k)}) - \Li^* u_{n(k)}}{\|-\G(\x_{n(k)}) - \Li^* u_{n(k)}\|}, \w_k = \frac{\Li^* u_{n(k)}}{\|\Li^* u_{n(k)}\|}. 
\end{equation}
Form the Cauchy–Schwarz inequality,
\begin{equation}
  \label{eq:low}
  (\forall k \in \N) \ \product[]{\vv_k}{\w_k} \geq - \|\vv_k\|\|\w_k\| = -1.
\end{equation}
In addition, we have
\begin{align}
  \product[]{\vv_k}{\w_k} 
  &= \product{\frac{-\G(\x_{n(k)}) - \Li^* u_{n(k)}}{\|-\G(\x_{n(k)}) - \Li^* u_{n(k)}\|}}{\frac{\Li^* u_{n(k)}}{\|\Li^* u_{n(k)}\|}} \\
  &= -\frac{\product{\G(\x_{n(k)})}{\Li^* u_{n(k)}}}{\|-\G(\x_{n(k)}) - \Li^* u_{n(k)}\|\|\Li^* u_{n(k)}\|} - \frac{\|\Li^* u_{n(k)}\|}{\|-\G(\x_{n(k)}) - \Li^* u_{n(k)}\|} \\
  &\leq \frac{\|\G(\x_{n(k)})\|}{\|-\G(\x_{n(k)}) - \Li^* u_{n(k)}\|} - \frac{\|\Li^* u_{n(k)}\|}{\|-\G(\x_{n(k)}) - \Li^* u_{n(k)}\|} \\
  &\leq \frac{M_{\G}}{\|-\G(\x_{n(k)}) - \Li^* u_{n(k)}\|} - \frac{\|\Li^* u_{n(k)}\|}{\|\G(\x_{n(k)})\| + \|\Li^* u_{n(k)}\|} \\
  \label{eq:prod}
  &\leq \frac{M_{\G}}{\|-\G(\x_{n(k)}) - \Li^* u_{n(k)}\|} - \frac{\|\Li^* u_{n(k)}\|}{M_{\G} + \|\Li^* u_{n(k)}\|}. 
\end{align}
By $\|-\G(\x_{n(k)}) - \Li^* u_{n(k)}\| \to \infty$ and $\|\Li^* u_{n(k)}\| \to \infty \ (k \to \infty)$, we have \linebreak$\lim_{k\to\infty} \frac{M_{\G}}{\|-\G(\x_{n(k)}) - \Li^* u_{n(k)}\|} - \frac{\|\Li^* u_{n(k)}\|}{M_{\G} + \|\Li^* u_{n(k)}\|} = -1$. 
Then, from \eqref{eq:low} and \eqref{eq:prod}, we have 
\begin{equation}
  \label{eq:inverse_dir}
  \lim_{k\to\infty}\product[]{\vv_k}{\w_k} = -1.
\end{equation} 

By Fact \ref{fact:subdifferential}(iv) with $\mathrm{ri}(\dom(\iota_D) - \mathrm{ran}(\Li)) = \mathrm{ri}(D - \mathcal{G}) = \operatorname{ri}(\mathcal{G}) = \mathcal{G}\ni 0$, we have 
\begin{equation}
  \label{eq:chain}
  \bigl(\partial (\iota_{\Li^{-1}(D)}) = \bigr)\partial (\iota_D \circ \Li) = \Li^* \circ \partial \iota_D \circ \Li. 
\end{equation}
Then, we have the following relation with $\A$ in \eqref{eq:A} and $\B$ in \eqref{eq:B}:
\begin{align}
  (\x, u) \in \Fix(\Talpha) 
  &\overset{\mathrm{Prop. \ \ref{prop:FBF_fix}(i)}}{\Leftrightarrow} 
  (\x, u) \in \zer(\A+\B) 
  \Leftrightarrow 
  \begin{cases}
    - \G(\x) - \Li^*u \in \partial \iota_{\C}(\x) \\
    \Li\x \in \partial \iota_D^* (u)
  \end{cases}\\
  &\overset{\mathrm{Fact \ \ref{fact:subdifferential}(i)}}{\Leftrightarrow}
  \begin{cases}
    - \G(\x) - \Li^*u \in \partial \iota_{\C}(\x)  \\
    u \in \partial \iota_D(\Li \x)
  \end{cases}
  \overset{\eqref{eq:chain}}{\Leftrightarrow}
  \begin{cases}
    - \G(\x) - \Li^*u \in \partial \iota_{\C}(\x)  \\
    \Li^*u \in \partial \iota_{\Li^{-1}(D)}(\x), 
  \end{cases}
\end{align}
where the last relation holds by \ref{item:bound_cond3}. 
By $(\x_{n(k)}, u_{n(k)})_{k\in\N} \subset \Fix(\Talpha)$, we have $- \G(\x_{n(k)}) - \Li^*u_{n(k)} \in \partial \iota_{\C}(\x_{n(k)})$ and $\Li^*u_{n(k)} \in \partial \iota_{\Li^{-1}(D)}(\x_{n(k)})$ for every $k \in \N$. 
By Fact \ref{fact:subdifferential}(iii),
\begin{equation}
  \label{eq:inclusion}
  (\forall k \in \N) \ 
  \vv_k \in \partial \iota_{\C}(\x_{n(k)}) \ \mathrm{and} \ 
  \w_k \in \partial\iota_{\Li^{-1}(D)}(\x_{n(k)}).
\end{equation}
Since $(\x_{n(k)})_{k\in\N}, (\vv_{k})_{k\in\N}$, and $(\w_{k})_{k\in\N}$ are bounded, there exist subsequences $(\x_{n(k(l))})_{l\in\N}$, $(\vv_{k(l)})_{l\in\N}$, and $(\w_{k(l)})_{l\in\N}$ such that $(\x_{{n(k(l))}}$, $\vv_{k(l)}, \w_{k(l)}) \to (\bar{\x}, \bar{\vv}, \bar{\w}) \in \bm{\mathcal{H}}^3 \ (l \to \infty)$. 
By \eqref{eq:inclusion} and \cite[Prop. 20.38(iii)]{BC2017}, we obtain $\bar{\vv} \in \partial \iota_{\C}(\bar{\x})$ and $\bar{\w} \in \partial \iota_{\Li^{-1}(D)}(\bar{\x})$. 
In addition, we have $\bar{\vv} \neq \bm{0}_{\bm{\mathcal{H}}}$ and $\bar{\w} \neq \bm{0}_{\bm{\mathcal{H}}}$ by $\|\bar{\vv}\| = \|\bar{\w}\| = 1$. 
From $\product{\bar{\vv}}{\bar{\w}} = \product{\lim_{l\to\infty} \vv_{k(l)}}{\lim_{l\to\infty} \w_{k(l)}}$\linebreak$= \lim_{l\to\infty}\product{\vv_{k(l)}}{\w_{k(l)}} = -1$ by \eqref{eq:inverse_dir}, we have $\bar{\vv} = -\bar{\w}(\Leftrightarrow \bar{\vv} + \bar{\w} = \bm{0}_{\bm{\mathcal{H}}})$, which contradicts the condition \ref{item:bound_cond4}. 
Therefore, $\Fix(\Talpha)$ is bounded. 
\qed

\section{Proof of Proposition \ref{prop:cond_sumrule_R-BC}}
For every $x \in S_1 \cap S_2$, $0_{\mathcal{H}} \in \partial \iota_{S_1}(x)$ and $0_{\mathcal{H}} \in \partial \iota_{S_2}(x)$ (see Fact \ref{fact:subdifferential}(iii)), and thus $\partial \iota_{S_1}(x) \cap - \partial \iota_{S_2}(x) \ni 0_{\mathcal{H}} \ (\forall x \in S_1 \cap S_2)$. 
Take $v \in \partial \iota_{S_1}(x) \cap - \partial \iota_{S_2}(x)$. 
In the following, we show $v = 0_{\mathcal{H}}$. 
By Fact \ref{fact:subdifferential}(iii), we have 
\begin{equation}
  \label{eq:bounded_1}
  (\forall y_1 \in S_1) \ \product[]{y_1 - x}{v} \leq 0 
\end{equation}
and 
\begin{align}
  \label{eq:bounded_2}
  (\forall y_2 \in S_2) \ \product[]{y_2 - x}{-v} \leq 0 \bigl( \Leftrightarrow (\forall y_2 \in S_2) \ \product[]{y_2 - x}{v} \geq 0 \bigr). 
\end{align}
From \eqref{eq:bounded_1} and \eqref{eq:bounded_2}, we have
\vspace{-0.5em}
\begin{equation}
  \label{eq:bounded_3}
  (\forall y \in S_1 \cap S_2) \ \product[]{y - x}{v} = 0. 
\end{equation}
\noindent\underline{Proof of \ref{item:QC_BC_1} $\Rightarrow$ \eqref{eq:QC_R}}. 
The condition \ref{item:QC_BC_1} ensures the existence of some $z \in \bigl(\operatorname{int} (S_1)\bigr) \cap S_2 (\subset S_1 \cap S_2)$. 
By $z \in \operatorname{int} (S_1)$, there exists $\delta > 0$ such that $B(z;\delta) \bigl( = z + B(0;\delta) \bigr) \subset S_1$. 
Since $z + w \in S_1$ and $z - w \in S_1$ hold for every $w \in B(0;\delta)$, we have 
\begin{equation}\hspace{-1em}
  \label{eq:bounded_4}
  \begin{aligned}
    &\bigl(\forall w \in B(0;\delta)\bigr)
    \begin{cases}
      \product[]{(z + w) - x}{v} \leq 0 \\
      \product[]{(z - w) - x}{v} \leq 0
    \end{cases}
    \Leftrightarrow
    \bigl(\forall w \in B(0;\delta)\bigr)
    \begin{cases}
      \product[]{z - x}{v} + \product[]{w}{v} \leq 0 \\
      \product[]{z - x}{v} - \product[]{w}{v}\leq 0
    \end{cases} \\
    &\overset{\eqref{eq:bounded_3}}{\Leftrightarrow}
    \bigl(\forall w \in B(0;\delta)\bigr)
    \begin{cases}
      \product[]{w}{v} \leq 0 \\
      \product[]{w}{v} \geq 0
    \end{cases}
    \Leftrightarrow
    \bigl(\forall w \in B(0;\delta)\bigr) \ \product[]{w}{v} = 0 
    \Leftrightarrow v = 0_{\mathcal{H}}. 
  \end{aligned}
\end{equation}
\underline{Proof of \ref{item:QC_BC_2} $\Rightarrow$ \eqref{eq:QC_R}}. 
The condition \ref{subitem:QC_ri_C1} ensures the existence of $z_1 \in \operatorname{ri} (S_1) \cap S_2(\subset S_1 \cap S_2)$. 
By $z_1 \in \operatorname{ri} (S_1)$, there exists $\delta>0$ such that $B(z_1;\delta) \cap \operatorname{aff} (S_1) \subset S_1$. 
Then,~we~have
\begin{align}
  B(z_1;\delta) \cap \operatorname{aff} (S_1) 
  &= \bigl(B(0;\delta) + z_1\bigr) \cap \bigl(\operatorname{span}(S_1 - z_1) + z_1\bigr)\\
  &= z_1 + \bigl(B(0;\delta) \cap \operatorname{span}(S_1 - z_1)\bigr) \subset S_1. 
\end{align}
Consequently, we have $z_1 + w \in S_1, z_1 - w \in S_1$ for every $w \in B(0;\delta) \cap \operatorname{span}(S_1 - z_1)$. 
Passing through a similar discussion to \eqref{eq:bounded_4}, 
we have $\product[]{w}{v} = 0$ for every $w \in B(0;\delta) \cap \operatorname{span}(S_1 - z_1)$. 
Since $\product[]{\lambda w}{v} = 0$ holds for every $w \in B(0;\delta) \cap \operatorname{span}(S_1 - z_1)$ and $\lambda \in \R$, we have 
\vspace{-0.5em}
\begin{equation}
  \label{eq:bounded_5}
  v \in \bigl(\operatorname{span}(S_1 - z_1)\bigr)^{\perp} \coloneqq \{v \in\mathcal{H} \mid (\forall \tilde{w} \in \operatorname{span}(S_1 - z_1)) \ \product[]{\tilde{w}}{v} = 0\}. \vspace{-0.5em}
\end{equation}
By a similar discussion, $v \in \bigl(\operatorname{span}(S_2 - z_2)\bigr)^{\perp}$ holds for $z_2 \in S_1\cap\operatorname{ri} (S_2)$ from the condition \ref{subitem:QC_ri_C2}. 
Then, we have 
\vspace{-0.5em}
\begin{equation}
  \label{eq:bounded_6}
  v \in \bigl(\operatorname{span}(S_1 - z_1)\bigr)^{\perp} \cap \bigl(\operatorname{span}(S_2 - z_2)\bigr)^{\perp} = \bigl( \operatorname{span}(S_1 - z_1) + \operatorname{span}(S_2 - z_2) \bigr)^{\perp}. 
  \vspace{-0.5em}
\end{equation}
In addition, from the condition \ref{subitem:QC_ri_aff}, we have 
\begin{equation}
  \vspace{-0.5em}
  \begin{aligned}
    \operatorname{aff}(S_1) + \operatorname{aff}(S_2) = \mathcal{H}
    &\Leftrightarrow
    \operatorname{span}(S_1 - z_1) + z_1 + \operatorname{span}(S_2 - z_2) + z_2 = \mathcal{H}\\
    &\Leftrightarrow
    \operatorname{span}(S_1 - z_1) + \operatorname{span}(S_2 - z_2) = \mathcal{H}\\
    &\Leftrightarrow
    \bigl( \operatorname{span}(S_1 - z_1) + \operatorname{span}(S_2 - z_2) \bigr)^{\perp} = \{0_{\mathcal{H}}\}.
  \end{aligned}
  \vspace{-0.5em}
\end{equation}
By \eqref{eq:bounded_6}, we obtain $v = 0_{\mathcal{H}}$.
\qed